\documentclass[journal]{IEEEtran}
\ifCLASSINFOpdf
\else
\fi
\usepackage{graphicx}      
\usepackage{amssymb,amsmath}
\usepackage{dsfont}
\usepackage{color}
\usepackage{epstopdf}        
\usepackage{amsthm}
\usepackage{amsfonts}
\usepackage{enumerate}
\usepackage{lscape}
\usepackage{multicol}
\usepackage{multirow}
\usepackage{calligra}
\usepackage{booktabs}
\usepackage{mathtools}
\usepackage{framed} 

\theoremstyle{remark}
\newtheorem{thm}{Theorem}
\newtheorem{cor}{Corollary}
\newtheorem{lemma}{Lemma}

\newtheorem{assumption}{Assumption}
\newtheorem{remark}{Remark}

\usepackage{graphics} 
\usepackage{epsfig} 
\usepackage{times} 
\usepackage{amsmath} 
\usepackage{amssymb}  
\usepackage{amsfonts}  
\usepackage{subfig}
\usepackage{epstopdf}
\usepackage{enumerate}
\usepackage{graphicx} 
\usepackage{color}
\usepackage{xcolor}
\usepackage{mathtools}
\usepackage{comment}

\newcommand*{\QEDB}{\hfill\ensuremath{\square}}%

\begin{document}
\title{Fixed-Time Extremum Seeking}


\author{Jorge I. Poveda, \textsl{Member, IEEE}, and Miroslav Krsti\'{c}, \textsl{Fellow, IEEE} 
\thanks{J. I. Poveda is with the Department of Electrical, Computer and Energy Engineering, University of Colorado, Boulder, CO, USA. E-mail: {\tt jorge.poveda@colorado.edu}.}
\thanks{ M. Krsti\'{c} is with the Department of Mechanical and Aerospace Engineering, University of California, San Diego, La Jolla, CA, USA. E-mail:   
        {\tt krstic@ucsd.edu}.}
\thanks{Research supported in part by CU Boulder - Autonomous Systems IRT Seed Grant No. 11005946, and NSF grants 1823983 $\&$ 1711373.   }
\thanks{The material in this paper was partially presented in the conference submissions \cite{ACCpovedaKrstic} and \cite{PovedaKrsticIFACWC20}.}
%
}


\maketitle

\begin{abstract}
We introduce a new class of extremum seeking controllers able to achieve fixed time convergence to the solution of optimization problems defined by static and dynamical systems. Unlike existing approaches in the literature, the convergence time of the proposed algorithms does not depend on the initial conditions and it can be prescribed a priori by tuning the parameters of the controller. Specifically, our first contribution is a novel gradient-based extremum seeking algorithm for cost functions that satisfy the Polyak-Lojasiewicz (PL) inequality with some coefficient $\kappa>0$, and for which the extremum seeking controller guarantees a fixed upper bound on the convergence time that is independent of the initial conditions but dependent on the coefficient $\kappa $. Second, in order to remove the dependence on $\kappa$, we introduce a novel Newton-based extremum seeking algorithm that guarantees a fully assignable fixed upper bound on the convergence time, thus paralleling existing asymptotic results in Newton-based extremum seeking where the rate of convergence is fully assignable. Finally, we study the problem of optimizing dynamical systems, where the cost function corresponds to the steady state input-to-output map of a stable but unknown dynamical system. In this case, after a time scale transformation is performed, the proposed extremum seeking controllers achieve the same fixed upper bound on the convergence time as in the static case. Our results exploit recent gradient flow structures proposed by Garg and Panagou in \cite{fixed_time}, and are established by using averaging theory and singular perturbation theory for dynamical systems that are not necessarily Lipschitz continuous. We confirm the validity of our results via numerical simulations that illustrate the key advantages of the extremum seeking controllers presented in this paper. 
\end{abstract}
\vspace{-0.4cm}
\section{Introduction}
\label{Sec:introduction}
%
%
\PARstart{I}{n} several applications it is of interest to recursively minimize a particular cost function whose mathematical form is unknown and which is only accessible via measurements. For these types of problems, extremum seeking control (ESC) has shown to be a powerful technique with provable stability, convergence, and robustness guarantees \cite{KrsticBookESC,DerivativesESC,Moase,Durr:Lie}. The main underlaying idea behind extremum seeking control is to induce multiple time scales in the dynamics of the closed-loop system in oder to emulate the behavior of a target nominal optimization algorithm chosen a priori to solve a particular type of optimization problem. Typical target nominal algorithms include gradient descent for convex optimization \cite{tan06Auto,GrushkovskayaLie}, Newton-like methods \cite{PowerES}, Riemannian gradient flows for constrained optimization \cite{Poveda:15}, gradient descent with momentum \cite{HeavyBollES}, accelerated gradient descent \cite{zero_order_poveda_Lina}, hybrid gradient descent on manifolds \cite{Strizic:17_CDC}, gradient and Newton flows with delays \cite{ThiagoKrstic}, gradient descent with time-varying costs \cite{Grushkovskaya2017}, distributed gradient systems for multi-agent systems \cite{Suttner2017}, and switched gradient flows modeled as hybrid systems \cite{Poveda:16}. Other ES approaches based on parameter estimation in adaptive control have also been considered in \cite{Guay:03,Guay:15,AttaGuay18,PhasorGuay}. Modifications to improve accuracy in ES were also presented in \cite{improvedES}. The previous approaches have been successfully used in several engineering applications such as wind turbine control and power converter optimization \cite{PowerES}, resource allocation problems \cite{Poveda:15}, optimization of robotic systems \cite{Fish_ESC}, price seeking in dynamic markets \cite{Frihauf:12}, dynamic tolling for transportation systems \cite{PovedaCDC17_a}, and traffic light control \cite{Kutadinata:14_Traffic}, to just name a few.  However, even though significant progress has been made in the analysis and design of ESC during the last years, obtaining a good transient performance characterized by fast rates of convergence remains a persistent challenge. This issue has recently motivated the development of fast ESCs based on Netwon-like flows \cite{PowerES}, gradient algorithms with time-invariant momentum \cite{HeavyBollES}, discontinuous gradient descent with finite-time stability properties \cite[Sec. 6.1]{Poveda:16}, and hybrid gradient flows with momentum resetting \cite{zero_order_poveda_Lina}. An early use of non-smooth feedback in extremum seeking can be found in control laws (48) and (55) in \cite{NonC2Krstic}, as well as finite-time (initial condition-dependent) convergence in Figures 2, 4, and 6 of \cite{NonC2Krstic}. Nevertheless, while these recent approaches can significantly improve the transient performance of the system, the convergence properties of all these ESCs are still of (practical) \emph{asymptotic} nature in the sense that for each small neighborhood $\mathcal{N}_{\epsilon}$ that contains the set of optimizers, and for each compact set of initial conditions $K_0$, the controller can be tuned to guarantee convergence to $\mathcal{N}_{\epsilon}$ in finite time, but, in general, as $\mathcal{N}_{\varepsilon}$ shrinks and $K_0$ grows, the convergence time grows unbounded. 

On the other hand, recently there has been significant efforts in designing and analyzing  control, estimation, and optimization algorithms with non-asymptotic convergence properties. These algorithms guarantee convergence to the desired target in a \emph{fixed time} that is finite and independent of the initial conditions. Algorithms with fixed time stability properties have been recently studied in \cite{JaimeFixed_Time,finite_timeEngelTAC,Fixed_timeTAC,OldFinite_Time,PralyFinite_Time,fixed_time} and \cite{romeronips}. As shown in \cite{Fixed_timeTAC} and \cite{PralyFinite_Time}, this type of convergence property can be established via Lyapunov functions for a class of continuous-time dynamical systems, which has opened the door to novel optimization algorithms characterized by continuous vector fields with non-asymptotic convergence properties; e.g.,  \cite{fixed_time} and \cite{Garg_Inequalities}. However, to the best of our knowledge, all existing optimization algorithms with fixed-time stability properties are model-based and require access to the first or second derivatives of the cost functions. 

\vspace{-0.2cm}
\subsection{Contributions}
In this paper, we present a new class of extremum seeking dynamics able to achieve fixed-time (semi-global practical) convergence in static and dynamical systems. This type of convergence guarantees that a fixed upper bound on the convergence time can be prescribed a priori independently of the initial conditions. More precisely the contributions of this paper are the following:
\begin{enumerate}
\item We present a novel fixed-time gradient-based extremum seeking (FTGES) algorithm for cost functions described by smooth static maps that attain its minima and that satisfy the Polyak-Lojasiewicz (PL) inequality with some coefficient $\kappa>0$. For any function satisfying these properties, we establish the existence of tunable parameters that guarantee (semi-global practical) fixed time convergence to a neighborhood of the optimizers, with an upper bound on the convergence time that depends on the coefficient $\kappa$.  This establishes a clear advantage in comparison to existing accelerated and non accelerated gradient based extremum seeking controllers for which the convergence time grows unbounded with unbounded initial conditions of the optimizing state. 
\item In order to remove the dependence of the convergence time on the parameters of the cost function, we introduce a new fixed-time Newton-based extremum seeking algorithm (FTNES), for which the upper bound in the convergence time is completely assignable a priori. Given that Newton-based ESCs carry out the estimation of the inverse of the Hessian matrix of the cost function by using a Riccati differential equation that has multiple equilibria, our convergence results are local by nature. However, unlike existing results in the literature, the convergence time to an arbitrarily small neighborhood of the optimizer can be upper bounded by a positive number that is independent of the initial conditions of the optimizing state and the Hessian of the cost function, and which can be prescribed a priori by the designer. This exhibits a clear advantage in comparison to traditional Newton-based extremum seeking schemes whose convergence times depend on the initial conditions.
\item We extend our previous results to dynamic settings where the cost function corresponds to the steady-state input-to-output mapping of a stable dynamical system. In this case, we show that equivalent fixed-time convergence results can be obtained after a time scale transformation is performed. In turn, this implies that in the original time scale the fixed upper bound is scaled by the inverse of a gain of the extremum seeking controller that needs to be selected sufficiently small in order to guarantee stability of the closed-loop system. 
\item Finally, we show that, under a certain choice of parameters, the proposed algorithms reduce to the standard gradient-based and Newton-based extremum seeking dynamics, and we recover the existing asymptotic results found in the literature of ESC. In this way, the dynamics introduced in this paper can be seen as generalizations of the standard gradient-based and Newton-based ESC.
\end{enumerate}
In order to analyze the extremum seeking dynamics considered in this paper, we make use of averaging theorems developed for non-smooth and set-valued systems \cite{Wang:12_Automatica,averaging_singularHDS,zero_order_poveda_Lina}. This allows us to link the $\mathcal{K}\mathcal{L}$ bound of the average system with the $\mathcal{K}\mathcal{L}$ bound that characterizes the properties of the extremum seeking dynamics, thus establishing a semi-global practical fixed-time convergence property. We illustrate the performance of our ESCs via simulations in different scenarios of single-variable and multi-variable optimization problems, comparing the trajectories generated by the algorithm with the trajectories generated by the standard vanilla gradient-based and Newton-based ESCs of \cite{DerivativesESC} and \cite{Newton}, respectively. To the knowledge of the authors, our results correspond to the first averaging-based extremum seeking algorithms with fixed-time convergence properties for static maps and dynamical systems.

\subsection{Additional contributions with respect to the conference submissions \cite{ACCpovedaKrstic} and \cite{PovedaKrsticIFACWC20}}
In contrast to the preliminary results presented in the conference submissions \cite{ACCpovedaKrstic} and \cite{PovedaKrsticIFACWC20}, in this journal paper we present the complete stability proofs of the algorithms by exploiting continuity of the vector fields. We also present a detailed analysis of the existence of complete solutions, as well as novel closed-loop architectures, stability results, and numerical examples for fixed-time gradient-based and Newton-based extremum seeking controllers applied to \emph{dynamical systems}. We further establish connections with previous results in the literature, that illustrate the advantages and generality of the proposed algorithms. 
\subsection{Organization}
The rest of this paper is organized as follows: Section \ref{Sec2} introduces the notation and some preliminaries on dynamical systems. Sections \ref{Sec_Gradient} and \ref{Sec_Newton} present the fixed-time gradient-based and Newton-based extremum seeking algorithms, respectively, as well as their main convergence properties and stability analysis. Section \ref{Sec_Dynamic} presents the results of the extremum seeking algorithms applied to dynamical systems, and finally Section \ref{Sec_Conclusions} ends with the outlook and some conclusions. Our theoretical results are illustrated throughout the paper by means of numerical examples.
%
\section{Preliminaries}
\label{Sec2}
%
\subsection{Notation}
We denote by $\mathbb{R}$ ($\mathbb{R}_{>0}$) the set of real numbers (resp. positive real numbers), and by $\mathbb{Z}$ ($\mathbb{Z}_{>0}$) the set of integers (resp. positive integers). Given a compact set $\mathcal{A}\subset\mathbb{R}^n$ and a vector $z\in\mathbb{R}^n$, we use $|z|_{\mathcal{A}}:=\min_{s\in\mathcal{A}}\|z-s\|_2$ to denote the minimum distance of $z$ to $\mathcal{A}$.  We use $\mathbb{S}^1:=\{z\in\mathbb{R}^2:z^2_1+z_2^2=1\}$ to denote the unit circle in $\mathbb{R}^2$, and $r\mathbb{B}$ to denote a closed ball in the Euclidean space, of radius $r>0$, and centered at the origin.  A function $\alpha:\mathbb{R}_{\geq0}\to\mathbb{R}_{\geq0}$ is of class $\mathcal{K}_{\infty}$ if it is zero at zero, continuous, strictly increasing, and unbounded. A function $\beta:\mathbb{R}_{\geq0}\times\mathbb{R}_{\geq0}\to\mathbb{R}_{\geq0}$ is of class $\mathcal{K}\mathcal{L}$ if it is nondecreasing in its first argument, nonincreasing in its second argument, $\lim_{r\to0^+}\beta(r,s)=0$ for each $s\in\mathbb{R}_{\geq0}$, and  $\lim_{s\to\infty}\beta(r,s)=0$ for each $r\in\mathbb{R}_{\geq0}$. 
%
%
We define the matrix $\mathcal{D}\in\mathbb{R}^{n\times{2n}}$ as the binary matrix that maps a vector $z=[z_1,z_2,z_3,\ldots,z_{2n}]^\top\in\mathbb{R}^{2n}$ to a vector $\tilde{z}$ having only the odd components of $z$, i.e., $\tilde{z}=\mathcal{D}z=[z_1,z_3,z_5,\ldots,z_{2n-1}]^\top$.

\subsection{Dynamical Systems and Stability Notions}
In this paper, we consider constrained dynamical systems with state $x\in\mathbb{R}^n$ and dynamics of the form
\begin{equation}\label{ODE}
x\in C,~~~~\dot{x}=F(x),
\end{equation}
where $F:\mathbb{R}^n\to\mathbb{R}^n$ is a continuous function and $C\subset\mathbb{R}^n$ is a closed set. A solution to system \eqref{ODE} is a continuously differentiable function $x:\text{dom}(x)\to\mathbb{R}^n$ that satisfies: a) $x(0)\in C$; b) $x(t)\in C$ for all $t\in\text{dom}(x)$; and c) $\dot{x}(t)=f(x(t))$ for all $t\in\text{dom}(x)$. A solution is said to be complete if $\text{dom}(x)=[0,\infty)$. Given a compact set $\mathcal{A}\subset C$, system \eqref{ODE} is said to render $\mathcal{A}$ uniformly globally asymptotically stable (UGAS) if there exists a class $\mathcal{K}\mathcal{L}$ function $\beta$ such that every solution of \eqref{ODE} satisfies 
\begin{equation}\label{KL_bound}
|x(t)|_{\mathcal{A}}\leq \beta(|x(0)|_{\mathcal{A}},t),~~~\forall~t\in\text{dom}(x).
\end{equation}
%

In this paper we will also consider $\varepsilon$-perturbed or parameterized dynamical systems of the form
\begin{equation}\label{perturbed_ode}
x\in C_{\varepsilon},~~~\dot{x}=f_{\varepsilon}(x),
\end{equation}
whose stability properties can be established only in a semi-global practical way. In particular, a compact set $\mathcal{A}\subset C$ is said to be semi-globally practically asymptotically stable (SGPAS) as $\varepsilon\to0^+$ if there exists a class $\mathcal{K}\mathcal{L}$ function $\beta$ such that for each pair $\delta>\nu>0$ there exists $\varepsilon^*>0$ such that for all $\varepsilon\in(0,\varepsilon^*)$ every solution of \eqref{perturbed_ode} with $|x(0)|_{\mathcal{A}}\leq \delta$ satisfies 
\begin{equation}
|x(t)|_{\mathcal{A}}\leq \beta(|x(0)|_{\mathcal{A}},t)+\nu,~~~\forall~t\in\text{dom}(x).
\end{equation}
The notion of SGPAS can be extended to systems that depend on multiple parameters $\varepsilon=[\varepsilon_1,\varepsilon_2,\ldots,\varepsilon_{\ell}]^\top$. In this case, and with some abuse of notation, we say that the system \eqref{perturbed_ode} renders the set $\mathcal{A}$ SGPAS as $(\varepsilon_{\ell},\ldots,\varepsilon_2,\varepsilon_{1})\to0^+$, where the parameters are tuned in order starting from $\varepsilon_1$, i.e., the parameters $(\varepsilon_{\ell},\ldots,\varepsilon_2,\varepsilon_{1})$ may not be selected independently. This type of stability notion is standard in extremum seeking control, see \cite{DerivativesESC,Poveda:16}.

%
%

\subsection{Dynamic Oscillators}
Our model-free optimization algorithms make use of several sinusoid signals that facilitate the extraction of gradient and hessian-related information of the cost function. In order to model these excitation signals, we consider $n$ uncoupled linear dynamic oscillators evolving on the $n$-torus $\mathbb{S}^n=\mathbb{S}^1\times\mathbb{S}^1\times\ldots\times\mathbb{S}^1\subset\mathbb{R}^{2n}$ with overall state $\mu\in \mathbb{R}^{2n}$ and dynamics
\begin{equation}\label{oscillator}
\mu\in \mathbb{S}^n,~~~\dot{\mu}=-\frac{2\pi}{\varepsilon_1}\mathcal{R}_{\kappa}\mu,
\end{equation}
where $\varepsilon_1>0$. The matrix $\mathcal{R}_{\kappa}\in\mathbb{R}^{2n\times 2n}$ is block diagonal and parametrized by a vector of gains $\kappa=[\kappa_1,\kappa_2,\ldots,\kappa_n]^\top$. The $i^{th}$ diagonal block of $\mathcal{R}_{\kappa}$ is defined as
\begin{equation*}
\mathcal{R}_i:=\left[\begin{array}{cc}
0 & -\kappa_i\\
\kappa_i & 0
\end{array}\right],
\end{equation*}
where $\kappa_i>0$. Since the $n$ oscillators are uncoupled, the odd entries $\mu_i$ of the solutions $\mu$ generated by \eqref{oscillator} with $\mu(0)\in \mathbb{S}^n$ are given by
\begin{equation}\label{solutions_oscillator}
\mu_i(t)=\mu_{i}(0)\cos\left(\frac{2\pi}{\varepsilon_1} \kappa_it\right)+\mu_{i+1}(0)\sin\left(\frac{2\pi}{\varepsilon_1} \kappa_it\right),
\end{equation}
with $\mu_{i}(0)^2+\mu_{i+1}(0)^2=1$, for all $i\in\{1,3,5,\ldots,n-1\}$. Indeed, by the structure of the oscillators, the set $\mathbb{S}^n$ is forward invariant, i.e.,  if $\mu(0)\in \mathbb{S}^n$ then $\mu(t)\in\mathbb{S}^n$ for all $t\geq0$. Moreover, since no solution of \eqref{oscillator} is defined outside the $n$-torus, the set $\mathbb{S}^n$ is trivially globally attractive. Therefore, system \eqref{oscillator} actually renders UGAS the set $\mathbb{S}^n$. This property will facilitate the stability analysis of our algorithms via averaging theory for non-smooth systems \cite{Wang:12_Automatica}.

\section{Gradient-Based Fixed-Time Extremum Seeking for Static Maps}
\label{Sec_Gradient}
We start by considering the fixed-time extremum seeking problem for static maps using a gradient descent-based architecture. In particular, we consider the following unconstrained optimization problem
\begin{equation}\label{main_problem}
\min_{z\in\mathbb{R}^n}~~\phi(z),
\end{equation}
where $\phi:\mathbb{R}^n\to\mathbb{R}$ is an unknown mathematical function that satisfies $\inf_{z\in\mathbb{R}^n}\phi(z)>-\infty$. Our goal is to design a robust feedback-based optimization algorithm that steers the state $z$ to a neighborhood of the set of solutions of \eqref{main_problem} in a fixed time, by using only measurements of $\phi$. Since algorithms of this form have no access to the first or second derivatives of $\phi$, they are usually referred to as \emph{zero-order methods} or extremum seeking controllers. 
\subsection{Qualitative Properties of the Cost Functions}
In order to solve problem \eqref{main_problem}, we start by considering cost functions that satisfy the following assumption:
\begin{assumption}\label{assumption_1}
The function $\phi:\mathbb{R}^n\to\mathbb{R}$ is twice continuously differentiable, radially unbounded, has a unique minimizer $z^*\in\mathbb{R}^n$, and there exists a constant $\kappa>0$ such that the function $\phi$ satisfies the Polyak-Lojasiewicz (PL) inequality:
\begin{equation}\label{condition2}
\phi(z)-\phi^*\leq \frac{1}{2\kappa}|\nabla \phi(z)|^2,
\end{equation}
for all $z\in\mathbb{R}^n$, where $\phi^*:=\phi(z^*)$.  \QEDB
\end{assumption}
The PL inequality \eqref{condition2} is satisfied by any strongly convex function, i.e., any function $\phi$ that satisfies
\begin{equation*}
\phi(z_1)\geq \phi(z_2)+\nabla \phi(z_2)^\top(z_1-z_2)+\frac{\kappa}{2}|z_1-z_2|^2,
\end{equation*}
for all $z_1,z_2$. Twice continuously differentiable strongly convex functions also satisfy $\nabla^2 \phi(z)\geq \kappa I$, for all $z\in\mathbb{R}^n$, see \cite[Exercise 7.26]{Beck2014}.  Since under Assumption \ref{assumption_1} the function $\phi$ is radially unbounded, all the nonempty level sets of $\phi$ are compact.
%
%

\subsection{Gradient-Based Fixed-Time Dynamics}
\begin{figure}[t!]
  \centering
    \includegraphics[width=0.4\textwidth]{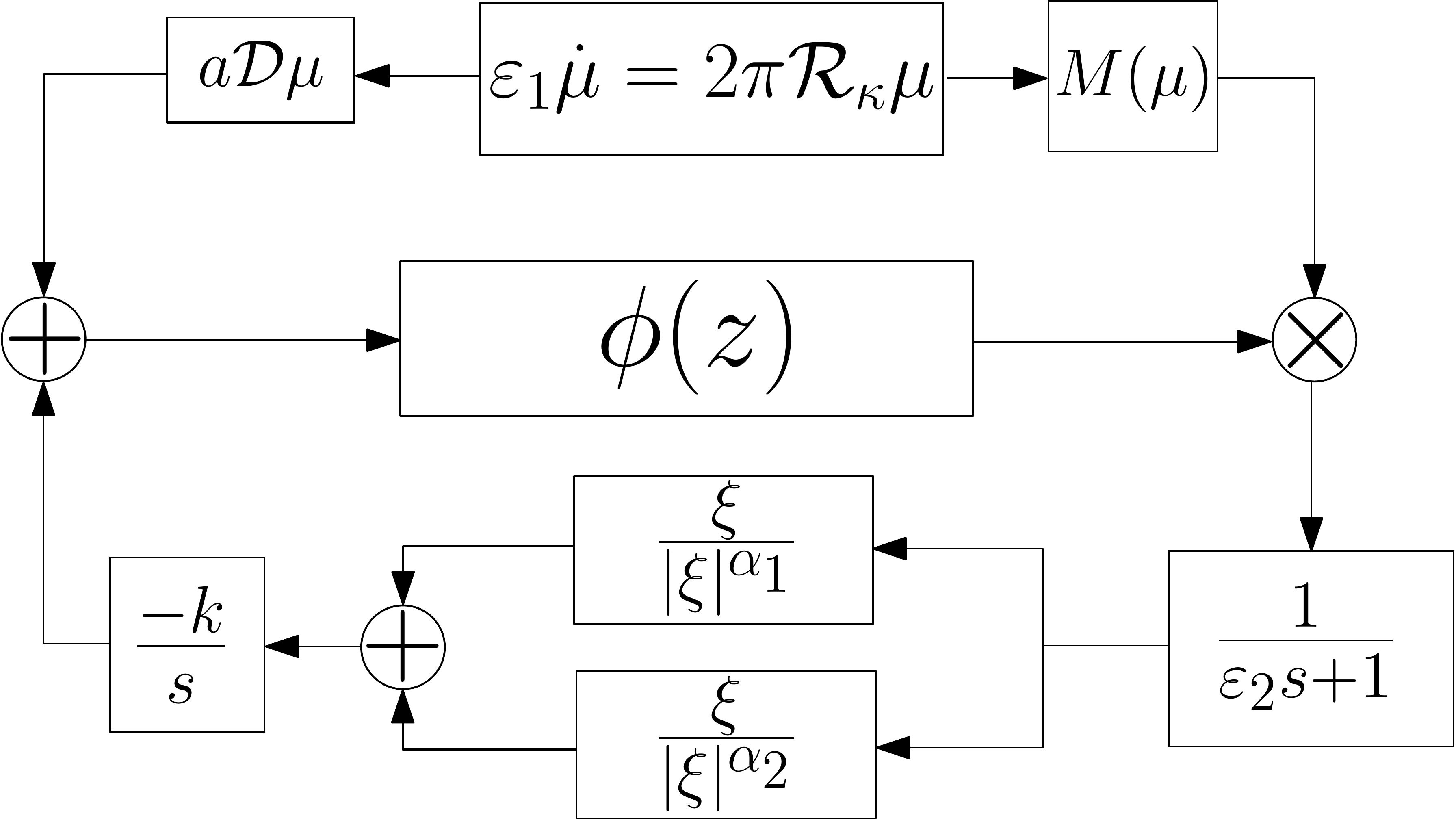}
    \caption{\small Scheme of the Fixed-Time Gradient-based Extremum Seeking (FTGES) algorithm for a static map $\phi$.}
    \label{figure11}
\end{figure}

In order to solve problem \eqref{main_problem} using only measurements of $\phi$, we consider  a Fixed-Time Gradient Extremum Seeking (FTGES) algorithm with state $(u,\xi,\mu)\in\mathbb{R}^n\times\mathbb{R}^n\times\mathbb{R}^{2n}$, evolving on the set 
\begin{equation}\label{flowset1}
(u,\xi,\mu)\in C:=\mathbb{R}^n\times \eta\mathbb{B} \times \mathbb{S}^n,
\end{equation}
with dynamics
\begin{align}\label{ES_dynamics1}
\left(\begin{array}{c}
\dot{u}\\\\
\dot{\xi}\\\\
\dot{\mu}
\end{array}\right)=-\left(\begin{array}{c}
k\left(\dfrac{\xi}{|\xi|^{\alpha_1}}+\dfrac{\xi}{|\xi|^{\alpha_2}}\right)\\
\dfrac{1}{\varepsilon_2}\Big(\xi-F_G(\phi,\mu)\Big)\\
\dfrac{2\pi}{\varepsilon_1}\mathcal{R}_{\kappa}\mu,
\end{array}\right),
\end{align}
where the right-hand side of $\dot{u}$ is defined as zero whenever $\xi=0$. The mapping $F_G$ is given by
\begin{equation}\label{mappingG}
F_{G}(\phi,\mu):=\phi(z)M(\mu),
\end{equation}
and where the argument of $\phi(\cdot)$, and the signal $M(\mu)$, are defined as
\begin{equation}\label{input}
z:=u+a\mathcal{D}\mu,~~~~~M(\mu):=\frac{2}{a}\mathcal{D}\mu,
\end{equation}
with $a\in\mathbb{R}_{\geq0}$ being a tunable parameter. The constants $\alpha_1$ and $\alpha_2$ are defined as
\begin{equation}\label{alphaconstants}
\alpha_1:=\frac{q_1-2}{q_1-1},~~~~\alpha_2:=\frac{q_2-2}{q_2-1},
\end{equation}
where $(q_1,q_2)\in\mathbb{R}_{>0}^2$ are tunable parameters that are said to be \emph{admissible} if they satisfy
\begin{equation*}
q_1\in(2,\infty)~~~\text{and}~~~~q_2\in(1,2).
\end{equation*}
In particular, admissible parameters $(q_1,q_2)$ guarantee that $\alpha_1$ is positive and $\alpha_2$ is negative. Figure \ref{figure11} shows a scheme illustrating the FTGES dynamics.

The proposed algorithm defines an extremum seeking controller with excitation signals $\mu$ generated by a linear oscillator, and a constrained low-pass filter with state $\xi$ evolving on a pre-defined compact set $\eta\mathbb{B}$. In practice, this compact set can be taken arbitrarily large by choosing a large $\eta$, and it is used only for the purpose of analysis in order to apply averaging results for nonsmooth singularly perturbed systems. In \eqref{ES_dynamics1}, the optimizing state $u$ is governed by learning dynamics that are designed to approximate a gradient flow with fixed-time convergence properties \cite{fixed_time}. However, said approximation requires suitable tuning of the parameters $(\varepsilon_1,\varepsilon_2, k, q_1,q_2, \kappa,a)$ that characterize the closed-loop system. The following assumption provides some general tuning guidelines.
\begin{assumption}\label{assumptionFTGES}
The parameters of the FTGES dynamics satisfy the following conditions:
\begin{enumerate}
\item The parameters $(q_1,q_2)$ are admissible.
\item The gain $k$ is positive, and $0<a\ll 1$.
\item For each $i\in\{1,2,\ldots,n\}$ the parameter $\kappa_i$ is a positive rational number, and $\kappa_i\neq \kappa_j$ for all $i\neq j$
\item The parameters $(\varepsilon_1,\varepsilon_2)$ satisfy  $0<\varepsilon_1\ll \varepsilon_2\ll 1$.
\end{enumerate}
\QEDB
\end{assumption}
Under Assumption \ref{assumptionFTGES}, the vector field describing the FTGES dynamics is continuous. In particular, since $1-\alpha_1>0$ and $1-\alpha_2>0$ we have that
\begin{align*}
\left|\lim_{\xi\to0}\frac{\xi}{|\xi|^{\alpha_1}}\right|&=\lim_{\xi\to0}\left|\frac{\xi}{|\xi|^{\alpha_1}}\right|
\leq \lim_{\xi\to0} |\xi|^{1-\alpha_1}=0.
\end{align*}
Similarly,
\begin{align*}
\left|\lim_{\xi\to0}\frac{\xi}{|\xi|^{\alpha_2}}\right|&=\lim_{\xi\to0}\left|\frac{\xi}{|\xi|^{\alpha_2}}\right|\leq \lim_{\xi\to0} |\xi|^{1-\alpha_2}=0.
\end{align*}
Therefore, the dynamics of the state $u$ are continuous at $\xi=0$. Continuity at points satisfying $\xi\neq0$ follows trivially by the structure of the right-hand side of \eqref{ES_dynamics1}. Existence of solutions follows then directly by \cite[Lemma 5.26]{HDS}.

The behavior of the FTGES can be roughly explained as follows: The dynamic oscillator generates sinusoid signals of the form \eqref{solutions_oscillator} with frequency proportional to $1/\varepsilon_1$. Thus, for values of $\varepsilon_1$ sufficiently small, the dynamics of $u$ and $\xi$ can be analyzed via averaging theory. By the construction of $F_G$ in \eqref{mappingG}, the average dynamics of $\xi$ will receive as input a perturbed estimation of the gradient of $\phi$, which is then feed to the dynamics of $u$ in order to carry out the optimization process. Provided $\varepsilon_2$ is sufficiently small, the qualitative behavior of the dynamics of $u$ can then be approximated by a simplified reduced system obtained via singular perturbation theory. Thus, under an appropriate tuning of the parameters $(\varepsilon_1,\varepsilon_2)$, the FTGES dynamics exhibit a multi-time scale behavior with three time scales, with the oscillator dynamics operating at the fastest time scale, and the optimization dynamics of $u$ operating at the slowest time scale. However, in contrast to existing smooth extremum seeking controllers, the FTGES dynamics cannot be analyzed using averaging theory and singular perturbation theory for Lipschitz continuous systems. Instead, we will use generalized averaging results developed for set-valued systems and hybrid dynamical systems with possibly non-unique solutions \cite{Wang:12_Automatica,zero_order_poveda_Lina}, which for ODEs only require continuity of the vector field.
\subsection{Main Result for the FTGES}
In order to state the main convergence result for the FTGES dynamics \eqref{ES_dynamics1}, we define for each admissible pair of parameters $(q_1,q_2)$ and each constant $\kappa$ satisfying the PL inequality \eqref{condition2}, the positive constants
\begin{equation}\label{gamma1}
\gamma_1=2^{\frac{8-3\alpha_1}{4}}\kappa^{\frac{2-\alpha_1}{2}},~~~\gamma_2=2^{\frac{8-3\alpha_2}{4}}\kappa^{\frac{2-\alpha_2}{2}},
\end{equation}
where $\alpha_1$ and $\alpha_2$ are defined in \eqref{alphaconstants}. Using \eqref{gamma1}, for each gain $k>0$, we define the value
\begin{equation}\label{fixed_time1}
T_G^*:=\frac{4}{k}\left(\frac{1}{\gamma_1\alpha_1}-\frac{1}{\gamma_2\alpha_2}\right)
\end{equation}
Since $(q_1,q_2)$ are admissible, the term inside the parenthesis is positive. Thus, for any desired $T_G^*>0$ equation \eqref{fixed_time1} can be satisfied by selecting the gain $k$ as
\begin{equation}\label{design_k}
k=\frac{4}{T_G^*}\left(\frac{1}{\gamma_1\alpha_1}-\frac{1}{\gamma_2\alpha_2}\right).
\end{equation}
We are now ready to state the first main result of this paper.
\vspace{0.1cm}

\begin{thm}\label{theorem1}
Suppose that $\phi$ satisfies Assumption \ref{assumption_1} and that Assumption \ref{assumptionFTGES} holds. Then, for any $T_G^*>0$ there exists admissible parameters $(q_1,q_2)$ and gain $k>0$ such that: for each pair $\delta>\nu>0$ there exists $\eta>0$ and $\varepsilon_2^*>0$ such that for each $\varepsilon_2\in(0,\varepsilon_2^*)$ there exists $a^*>0$ such that for each $a\in(0,a^*)$ there exists $\varepsilon_1^*>0$ such that for each $\varepsilon_1\in(0,\varepsilon_1^*)$ the FTGES dynamics \eqref{flowset1}-\eqref{ES_dynamics1} with $(u(0),\xi(0),\mu(0))\in \left(\{z^*\}+\delta\mathbb{B}\right)\times \eta\mathbb{B}\times\mathbb{S}^n$ generates solutions with unbounded time domain, and each of these solutions satisfies $|z(t)-z^*|\leq \nu$, for all $t\geq T_G^*$.  \QEDB
\end{thm}
%
\textsl{Proof:} In order to analyze the FTGES dynamics, we will use averaging tools for dynamical systems that are not necessarily Lipschitz continuous, e.g., \cite{Wang:12_Automatica,zero_order_poveda_Lina}. A key feature of these tools is that the $\mathcal{K}\mathcal{L}$ bound that characterizes the rate of convergence of the ``slow state'' in a singularly perturbed system is completely characterized by the $\mathcal{K}\mathcal{L}$ bound of its \emph{average system}. Based on this, we will divide the proof in two steps. First, by studying the stability properties of the compact set $\mathcal{A}\times \mathbb{S}^n$, where $\mathcal{A}:=\{z^*\}\times \eta\mathbb{B}$, we will establish a SGPAS result for the FTGES dynamics, which preserves the $\mathcal{K}\mathcal{L}$ bound of the average dynamics in the slowest time scale. After this, we will exploit the SGPAS result and the linearity of the low pass filter dynamics to show that the restriction of $\xi$ to the compact set $\eta\mathbb{B}$ is inconsequential for compact sets of initial conditions $L\mathbb{B}$ whenever $\eta\gg L$ is sufficiently large.

\vspace{0.1cm}
\textsl{Step 1:} We start with the following lemma, which follows directly by computation:
\begin{lemma}\label{properties_oscillator}
Under item 3) of Assumption \ref{assumptionFTGES} there exists a $T>0$ such that every solution $\mu$ of the oscillator \eqref{oscillator} with $\varepsilon_1=1$ satisfies
\begin{subequations}\label{orthogonality}
\begin{align}
&\frac{1}{\ell T}\int_{0}^{\ell T}\tilde{\mu}(t)\tilde{\mu}(t)dt=\frac{1}{2}I_n,\\
&\frac{1}{\ell T}\int_{0}^{\ell T}\tilde{\mu}(t)dt=0.
\end{align}
\end{subequations}
for all $\ell\in\mathbb{Z}_{\geq1}$, where $\tilde{\mu}=\mathcal{D}\mu$. \QEDB
\end{lemma}
By using the properties \eqref{orthogonality}, we can average the dynamics of $u$ and $\xi$ along the solutions $\mu$ of the oscillator.  To find the average dynamics of system \eqref{ES_dynamics1}, and for $a>0$ sufficiently small, consider a Taylor expansion of $\phi(u+a\tilde{\mu})$ around $u$: 
\begin{equation*}
\phi(u+a\tilde{\mu})=\phi(u)+a\tilde{\mu}^\top \nabla\phi(u)+O(a^2),
\end{equation*}
where we used $\tilde{\mu}:=\mathbb{D}\mu$ to shorten notation. Substituting in \eqref{mappingG}  and using the properties \eqref{orthogonality}, we obtain the following average dynamics with state $(u^a,\xi^a)$:
\begin{align}\label{ES_dynamics_average}
\left(\begin{array}{c}
\dot{u}^a\\\\
\dot{\xi}^a
\end{array}\right)=-\left(\begin{array}{c}
k\left(\dfrac{\xi^a}{|\xi^a|^{\alpha_1}}+\dfrac{\xi^a}{|\xi^a|^{\alpha_2}}\right)\\
\dfrac{1}{\varepsilon_2}\Big(\xi^a-\tilde{F}_G(\nabla \phi)\Big)
\end{array}\right)
\end{align}
which are defined in the set
\begin{equation*}
(u^a,\xi^a)\in \mathbb{R}^n\times \eta\mathbb{B}
\end{equation*}
and where the function $\tilde{F}_G(\nabla \phi)$ is given by
\begin{equation*}
\tilde{F}_{G}(\nabla \phi):=\nabla \phi(u^a)+O(a).
\end{equation*}
System \eqref{ES_dynamics_average} is an $O(a)$-perturbed version of a \emph{nominal average} system with dynamics
\begin{align}\label{ES_dynamics_average_nominal}
\left(\begin{array}{c}
\dot{u}^a\\\\
\dot{\xi}^a
\end{array}\right)=-\left(\begin{array}{c}
k\left(\dfrac{\xi^a}{|\xi^a|^{\alpha_1}}+\dfrac{\xi^a}{|\xi^a|^{\alpha_2}}\right)\\
\dfrac{1}{\varepsilon_2}\big(\xi^a-\nabla \phi(u^a)\big)
\end{array}\right).
\end{align}
Indeed, since system \eqref{ES_dynamics_average_nominal} is continuous, by \cite[Thm. 7.21]{HDS} suitable stability properties will be preserved in a semi-global practical way provided $a$ is sufficiently small.  Therefore, we proceed to analyze the stability and convergence properties of the nominal system \eqref{ES_dynamics_average_nominal}.

For $\varepsilon_2$ sufficiently small, system \eqref{ES_dynamics_average_nominal} is in singular perturbation form, with the dynamics of $\xi^a$ acting as fast dynamics. To find the boundary layer dynamics of this system, we introduce the new time scale $\tau=t/\varepsilon_2$, which leads to 
\begin{align*}
\left(\begin{array}{c}
\dfrac{\partial u^a}{\partial \tau}\vspace{0.1cm}\\
\dfrac{\partial \xi^a}{\partial \tau}
\end{array}\right)=-\left(\begin{array}{c}
\varepsilon_2k\left(\dfrac{\xi^a}{|\xi^a|^{\alpha_1}}+\dfrac{\xi^a}{|\xi^a|^{\alpha_2}}\right)\\\
\xi^a-\nabla \phi(u^a)
\end{array}\right).
\end{align*}
The boundary layer dynamics are obtained by setting $\varepsilon_2=0$, that is $\frac{\partial u_{bl}^a}{\partial \tau}=0$ and
\begin{align}\label{ES_dynamics_average_bl}
\xi_{bl}\in\eta\mathbb{B},~~~\dfrac{\partial \xi_{bl}^a}{\partial \tau}=-\Big(\xi_{bl}^a-\nabla \phi(u_{bl}^a)\Big).
\end{align}
For each fixed $u_{bl}$, i.e., $\dot{u}_{bl}=0$, the solutions of this system satisfy
\begin{equation*}
\big|\xi_{bl}^{a}(t)-\nabla\phi(u_{bl}^a)\big|\leq c_1 e^{-c_2 t}\big|\xi_{bl}^{a}(0)-\nabla\phi(u_{bl}^a)\big|,
\end{equation*}
for all $t\in\text{dom}(\xi_{bl},u_{bl})$, for some $c_1,c_2>0$. Thus, the singularly perturbed system \eqref{ES_dynamics_average_nominal} has a well-defined reduced system \cite[Ex. 1]{Wang:12_Automatica}. This reduced system is obtained by substituting $\xi^a$ by its steady-state value $\nabla \phi(u^a)$, which leads to the following unconstrained dynamics with state $u_r$:
\begin{align}\label{ES_dynamics_average_nominal_reduced}
u_r\in\mathbb{R}^n,~~~\dot{u}_r=-k\left(\dfrac{\nabla \phi(u_r)}{|\nabla \phi(u_r)|^{\alpha_1}}+\dfrac{\nabla \phi(u_r)}{|\nabla \phi(u_r)|^{\alpha_2}}\right).
\end{align}
In order to analyze system \eqref{ES_dynamics_average_nominal_reduced}, we can follow the ideas of \cite{fixed_time} by considering the Lyapunov function
\begin{equation*}
V_G(u_r)=\frac{1}{2}(\phi(u_r)-\phi^*)^2,
\end{equation*}
which, under Assumption 1, satisfies $V(u_r)>0$ for all $u_r\neq z^*$, and $V(u_r)=0$ if an only if $u_r=z^*$. Moreover, all the level sets of  $V(u_r)$ are bounded. The derivative of $V(u_r)$ along the solutions of \eqref{ES_dynamics_average_nominal_reduced} is given by
\begin{equation*}
\dot{V}_G(u_r)=-k(\phi(u_r)-\phi^*)\bigg(|\nabla \phi(u_r)|^{\tilde{\alpha}_1}+|\nabla \phi(u_r)|^{\tilde{\alpha}_2}\bigg),
\end{equation*}
where $\tilde{\alpha}_1=2-\alpha_1>0$ and $\tilde{\alpha}_2=2-\alpha_2>0$. Using the PL inequality \eqref{condition2} we have that $-|\nabla \phi(u_r)|\leq -\sqrt{2\kappa}\big(\phi(u_r)-\phi^*\big)^{\frac{1}{2}}$. Thus, the derivative $\dot{V}_G$ can be upper bounded as
\begin{align*}
\dot{V}_G(u_r)&\leq -k\Big((2\kappa)^{\frac{\tilde{\alpha}_1}{2}}(\phi(u_r)-\phi^*)^{\frac{2+\tilde{\alpha}_1}{2}}\\
&~~~~~~~~~~~~~~~+(2\kappa)^{\frac{\tilde{\alpha}_2}{2}}(\phi(u_r)-\phi^*)^{\frac{2+\tilde{\alpha}_2}{2}}\Big),\\
&=-k\Big( (2\kappa)^{\frac{\tilde{\alpha}_1}{2}}   (2V(u_r))^{\frac{2+\tilde{\alpha}_1}{4}}+(2\kappa)^{\frac{\tilde{\alpha}_2}{2}} (2V(u_r))^{\frac{2+\tilde{\alpha}_2}{4}} \Big),\\
&=-k \Big(c_1V_G(u_r)^{\gamma_1}+c_2V_G(u_r)^{\gamma_2} \Big)<0,
\end{align*}
for all $u_r\neq z^*$, where
\begin{align*}
c_1:=2^{\frac{2+3\tilde{\alpha_1}}{4}}\kappa^{\frac{\tilde{\alpha}_1}{2}}>0, \\
c_2:=2^{\frac{2+3\tilde{\alpha_2}}{4}}\kappa^{\frac{\tilde{\alpha}_2}{2}}>0,
\end{align*}
and 
\begin{align*}
\gamma_1:&=\frac{2+\tilde{\alpha}_1}{4} \in(0,1), \\
\gamma_2:&=\frac{2+\tilde{\alpha}_2}{4}>1.
\end{align*}
This establishes UGAS of $z^*$ for the reduced dynamics. Moreover, by \cite[Lemma 1]{Fixed_timeTAC}, the Lyapunov function evaluated along the solutions of \eqref{ES_dynamics_average_nominal_reduced} satisfies $V_G(u_r(t))=0$ for all $t\geq T_G^*$, where $T_G^*$ is given by \eqref{fixed_time1}. By the definition of $V_G(u_r)$ and Assumption \ref{assumption_1}, this implies that $u_r(t)=z^*$ for all $t\geq T_G^*$. 

The UGAS property implies the existence of a $\mathcal{K}\mathcal{L}$ bound $\beta_u$ such that all solutions of \eqref{ES_dynamics_average_nominal_reduced} satisfy
\begin{equation}\label{KLbound}
|u_r(t)-z^*|\leq \beta_u(|u_r(0)-z^*|,t),~~~~~\forall~~t\geq0.
\end{equation}
As noted in \cite[Thm. 3.40]{HDS}, the function $\beta_u$ can be constructed as $\beta_u(r,s)=\max\{0,\beta_0(r,s)\}$, where:
\begin{align}\label{KLexplicit}
\beta_{0}(r,s):=&\sup\Big\{|u_r(t)-z^*|:~u_r~\text{is a solution of \eqref{ES_dynamics_average_nominal_reduced}}\notag\\
 &~~~~~~\text{with}~u_r(0)=u_0,~|u_0-z^*|\leq r,~t\geq s\Big\}.
\end{align}
Since the solutions of \eqref{ES_dynamics_average_nominal_reduced} exists and are defined for all $t\geq0$, the supremum is well defined and the function is nondecreasing in $r$ and nonincreasing in $s$. Moreover, it satisfies $\lim_{r\to0^+}\beta_u(r,s)=0$ for each $s\in\mathbb{R}_{\geq0}$ since $\beta_u(r,s)\leq \alpha(r)$, where $\alpha\in\mathcal{K}_{\infty}$ is generated by the uniform global stability property. Finally, $\lim_{s\to\infty}\beta_u(r,s)=0$ for each $r$ due to the finite time convergence property, where the limit can be taken independent of $r$. Therefore, $\beta_u(r,s)$ is a valid class $\mathcal{K}\mathcal{L}$ function that bounds every solution of \eqref{ES_dynamics_average_nominal_reduced}, and by construction it satisfies
\begin{equation}\label{finite_bound}
\beta_u(|u_r(0)-z^*|,t)=0
\end{equation}
for all $t\geq T_G^*$ and all $u_r(0)\in\mathbb{R}^n$.

Having established UGAS and fixed time convergence for the reduced dynamics \eqref{ES_dynamics_average_nominal_reduced}, we can now use \cite[Thm. 2]{Wang:12_Automatica} to establish 
that system \eqref{ES_dynamics_average_nominal} renders the compact set $\mathcal{A}:=\{z^*\}\times \eta\mathbb{B}$ SGPAS as $\varepsilon_2\to 0^+$ with $\mathcal{K}\mathcal{L}$ bound $\beta_u$. Moreover, by the definition of solutions we have that $|\xi^a(t)|_{\eta\mathbb{B}}=0$ for all $t\in\text{dom}(u^a,\xi^a)$, which implies that $|\zeta^a(t)|_{\mathcal{A}}:=|u^a(t)-z^*|$, where $\zeta^a=(u^{a\top},\xi^{a\top})^\top$. Thus, for each $\delta>\nu>0$ there exists $\varepsilon_2^*>0$ such that for all $\varepsilon_2\in(0,\varepsilon^*_2)$ every solution of system \eqref{ES_dynamics_average_nominal} with $\zeta^a(0)\in (\{z^*\}+\delta\mathbb{B})\times \eta\mathbb{B}$ satisfies the following bound:
\begin{equation}\label{first_KLbound}
|\zeta^a(t)|_{\mathcal{A}}\leq \beta_u(|\zeta^a(0)|_{\mathcal{A}},t)+\nu,
\end{equation}
for all $t\in\text{dom}(\zeta^a)$. Having established the bound \eqref{first_KLbound}, we can now exploit the structural robustness properties of system \eqref{ES_dynamics_average_nominal}, which follow by the continuity of the right hand side of the dynamics, in order to establish via \cite[Thm. 7.21]{HDS} and \cite[Prop. 6]{zero_order_poveda_Lina}, that the $O(a)$-perturbed system \eqref{ES_dynamics_average} renders the same compact set $\mathcal{A}$ SGPAS as $(a,\varepsilon_2)\to 0^+$ with $\mathcal{K}\mathcal{L}$ bound $\beta_u$, i.e., for each pair $\delta>\nu>0$ there exists $\varepsilon_2^*>0$ such that for all $\varepsilon_2\in(0,\varepsilon_2^*)$ there exists $a^*>0$ such that for all $a\in(0,a^*)$ every solution of the perturbed system \eqref{ES_dynamics_average} with $u^a(0)\in \{z^*\}+\delta\mathbb{B}$ satisfies a bound of the form \eqref{first_KLbound}. 

Finally, since the fast oscillator dynamics of \eqref{ES_dynamics1} render UGAS the set $\mathbb{S}^n$ and generate a well-defined average system corresponding to \eqref{ES_dynamics_average}, by averaging results for perturbed non-Lipschitz systems \cite[Thm. 9]{zero_order_poveda_Lina} we obtain that the FTGES dynamics render the set $\mathcal{A}\times\mathbb{S}^n$ SGPAS as $(\varepsilon_1,a,\varepsilon_2)\to 0^+$ with $\mathcal{K}\mathcal{L}$ bound $\beta_u$. In particular, this establishes that for each $k>0$, each tuple of admissible parameters $(q_1,q_2)$, and each pair $\delta>\nu>0$ there exists $\varepsilon_2^*>0$ such that for each $\varepsilon_2\in(0,\varepsilon_2^*)$ there exists $a^*\in(0,\nu/2)$ such that for each $a\in(0,a^*)$ there exists $\varepsilon_1^*>0$ such that for each $\varepsilon_1\in(0,\varepsilon_1^*)$ each solution of the FTGES dynamics \eqref{flowset1}-\eqref{ES_dynamics1} satisfies the bound
\begin{equation}\label{KLLL}
|\zeta(t)|_{\mathcal{A}}\leq \beta_u(|\zeta(t)|_{\mathcal{A}},t)+\frac{\nu}{2},
\end{equation}
for all $t\in\text{dom}(\zeta,\mu)$, where $\zeta:=(u^\top,\xi^\top)^\top$. Given that by definition of solutions we have that $|\xi|_{\eta\mathbb{B}}=0$ and therefore $|\zeta|_{\mathcal{A}}=|u-z^*|$, it then follows that each solution of the FTGES dynamics \eqref{flowset1}-\eqref{ES_dynamics1} satisfies the bound
\begin{equation}\label{KLL}
|u(t)-z^*|\leq \beta_u(|u(t)-z^*|,t)+\frac{\nu}{2},
\end{equation}
for all $t\in\text{dom}(\zeta,\mu)$, where $\zeta:=(u^\top,\xi^\top)^\top$. Since the $\mathcal{K}\mathcal{L}$ bound $\beta_u$ satisfies the property of equation \eqref{finite_bound}, we obtain that $|u(t)-z^*|\leq\nu/2$ for all $t\in\text{dom}(u,\xi,\mu)$ such that $t\geq T_G^*$, namely, solutions with a time domain larger than $T^*_{G}$ achieve fixed-time convergence to a $\nu$-neighborhood of $z^*$. Using \eqref{input}, the triangle inequality, and the fact that $a\in(0,\nu/2)$, we obtain that $|z(t)-z^*|\leq\nu$ for all $t\in\text{dom}(u,\xi,\mu)$ such that $t\geq T_G^*$.

\vspace{0.1cm}
\textsl{Step 2:} We now exploit the $\mathcal{K}\mathcal{L}$ bound $\beta_u$ of system \eqref{ES_dynamics_average_nominal_reduced} in order to establish the existence of complete solutions for the FTGES dynamics from arbitrarily large compact sets $\left(\{z^*\}+\delta\mathbb{B}\right)\times L\mathbb{B}\times\mathbb{S}^n$ of initial conditions. Without loss of generality we assume that $\nu<1$. Due to the bound \eqref{KLL}, the fact that for any $\eta>0$ the set $\eta\mathbb{B}$ is compact, and the forward invariance of the compact set $\mathbb{S}^n$, by the results of Step 1 the FTGES dynamics have no finite escape times from compact sets of initial conditions $\{z^*\}+\delta\mathbb{B}$ and for suitable choices of parameters $(\varepsilon_2,a,\varepsilon_1)$. Thus, any solution of \eqref{flowset1}-\eqref{ES_dynamics1} with $\text{dom}(u,\xi,\mu)<\infty$ must stop due to $\xi$ leaving the set $\eta\mathbb{B}$. In order to establish the existence of solutions that do not stop, we note that for any uniformly bounded input $s$ the solutions of the linear dynamics $\dot{\xi^a}=\frac{1}{\varepsilon_2}\left(-\xi^a+s(t)\right)$, satisfy the bound
\begin{equation}\label{linear_bound}
|\xi^a(t)|\leq \exp\left(-\frac{1}{\varepsilon_2}t\right)|\xi^a(0)|+\sup_{t\geq0} |s(t)|,
\end{equation}
for all $t\geq0$ and all $\varepsilon_2>0$, see \cite[pp. 174]{khalil}. Fix the admissible parameters $(k,q_1,q_2)$ which define the fixed time $T_G^*$. Fixed the constants $\delta>\nu$ with $\nu\in(0,1)$, and define the set 

\vspace{-0.4cm}
\begin{small}
\begin{equation*}
\tilde{K}:=\left\{u\in\mathbb{R}^n:|u-z^*|\leq \beta_u\left(\max_{u_0\in \{z^*\}+\delta\mathbb{B}}|u_0-z^*|,0\right)+1\right\}.
\end{equation*}
\end{small}\noindent
By construction, this set is compact since without loss of generality $\beta_u$ can be taken to be continuous \cite[pp. 69]{HDS}. Thus, there exists $M>0$ such that $\tilde{K}\subset M\mathbb{B}$. Since $\phi$ is $C^2$, we have that $|u|\leq M$ implies $|\nabla \phi(u)|\leq M'$ for some $M'>0$. Let $L>0$ and consider the compact set $L\mathbb{B}\subset\mathbb{R}^n$. Let $M''>L$ and define the constant $\Gamma:= M''+M'$. Let $\eta:=2\Gamma$  and let Step 1 generate the values of $(\varepsilon^*_1,a^*,\varepsilon_2^*)$ that induce the bound \eqref{KLL}. Then, by the results of Step 1, for all $\varepsilon_2\in(0,\varepsilon_2^*)$ every solution of the singularly-perturbed system \eqref{ES_dynamics_average_nominal} with $u^a(0)\in \{z^*\}+\delta\mathbb{B}$ and $\xi^a(0)\in L\mathbb{B}$ generates trajectories $u^a$ that satisfy a bound of the form $|u^a(t)-z^*|\leq \beta_u(|u^a(0)-z^*|,0)+\nu/2$. Moreover, by linearity of the dynamics of the low-pass filter, the trajectories $\xi^a$ satisfy a bound of the form $\eqref{linear_bound}$, which, by the construction of $\Gamma\mathbb{B}$ and the choice of $\eta$,  implies that $\xi^a(t)\in \Gamma\mathbb{B}\subset \text{int}(\eta\mathbb{B})$ for all $t\geq0$. Thus, said solutions satisfy $\text{dom}(u^a,\xi^a)=[0,\infty)$. Since $\nabla \phi$ is locally Lipschitz, $\nabla \phi(z^*)=0$, and since $|u^a(t)-z^*|\leq \nu$ for all $t\geq T_G^*$, it follows that $|\nabla\phi(u^a(t))|\leq L_{z^*}\nu$ for all $t\geq T_G^*$, where $L_{z^*}>0$. Thus, by using the bound \eqref{linear_bound} with $\sup_{t\geq0} |s(t)|=L_{z^*}\nu$ it follows that trajectories $\xi^a$ of \eqref{ES_dynamics_average_nominal} with $\xi^a(0)\in L\mathbb{B}$ converge to a $L_{z^*}\nu$-neighborhood of zero.

Finally, we use $\epsilon$-closeness of solutions on compact time domains (\cite[Def. 5.23]{HDS}) between perturbed and unperturbed ODEs with a continuous right-hand side in order to establish the existence of complete solutions for the original FTGES dynamics: By \cite[Prop. 6.34]{HDS}, there exists $a^{**}>0$ such that for all $a\in(0,\min\{a^*,a^{**}\})$ and each solution of the $O(a)$-perturbed average dynamics \eqref{ES_dynamics_average} there exists a solution of the nominal average dynamics \eqref{ES_dynamics_average_nominal} that is $\epsilon$-close on compact time domains. Thus, since $\text{int}(\eta\mathbb{B})$ is an open set there exists $\epsilon>0$ such that for any $\bar{T}>T_G^*$ there exists $a^{***}>0$ sufficiently small such that for all $a\in(0,a^{***})$ the trajectories $\xi^a$ of system \eqref{ES_dynamics_average} with $(u^a(0),\xi^a(0))\in \left(\{z^*\}+\delta\mathbb{B}\right)\times L\mathbb{B}$ also satisfy $\xi^a(t)\in (\Gamma+\epsilon/2)\mathbb{B}\subset \text{int}(\eta\mathbb{B})$ for all $t\in[0,\bar{T}]$. By using again closeness of solutions between systems \eqref{ES_dynamics1} and \eqref{ES_dynamics_average}, we obtain the existence of $\varepsilon_1^{**}>0$ such that for all $\varepsilon_1\in(0,\min\{\varepsilon_1^*,\varepsilon_1^{**}\})$ for each trajectory of the original FTGES dynamics with $(u(0),\xi(0),\mu(0))\in\left(\{z^*\}+\delta\mathbb{B}\right)\times L\mathbb{B}\times\mathbb{S}^n$ there exists a solution of the average dynamics \eqref{ES_dynamics_average} that is $\epsilon/2$-close on compact time domains. Thus, the trajectories $(u,\xi)$ of the FTGES dynamics \eqref{ES_dynamics1} with $(u(0),\xi(0))\in \left(\{z^*\}+\delta\mathbb{B}\right)\times L\mathbb{B} $ also satisfy  $\xi(t)\in (\Gamma+\epsilon)\mathbb{B}\subset \text{int}(\eta\mathbb{B})$ for all $t\in[0,\bar{T}]$, where $\bar{T}\geq T^*_G$. Convergence of $\xi$ to an $L_{z^*}\nu$-neighborhood of the origin follows now directly by the bounds \eqref{linear_bound} and \eqref{KLL}, and the $\epsilon$-closeness of solutions on compact time domains between solutions of the FTGES dynamics and system \eqref{ES_dynamics_average_nominal}. In turn, this implies that $\xi(t)\in \text{int}(\eta\mathbb{B})$ for all $t\geq 0$, which establishes that every solution of the FTGES dynamics from the compact set of initial conditions $(\{z^*\}+\delta\mathbb{B})\times L\mathbb{B}\times\mathbb{S}^n\subset \mathbb{R}^n\times \eta\mathbb{B}\times\mathbb{S}^n$ satisfies $\text{dom}(u,\xi,\mu)=[0,\infty)$. \null\hfill\null $\blacksquare$

\subsection{Discussion and Numerical Examples}
\begin{figure}[t!]
  \centering
    \includegraphics[width=0.5\textwidth]{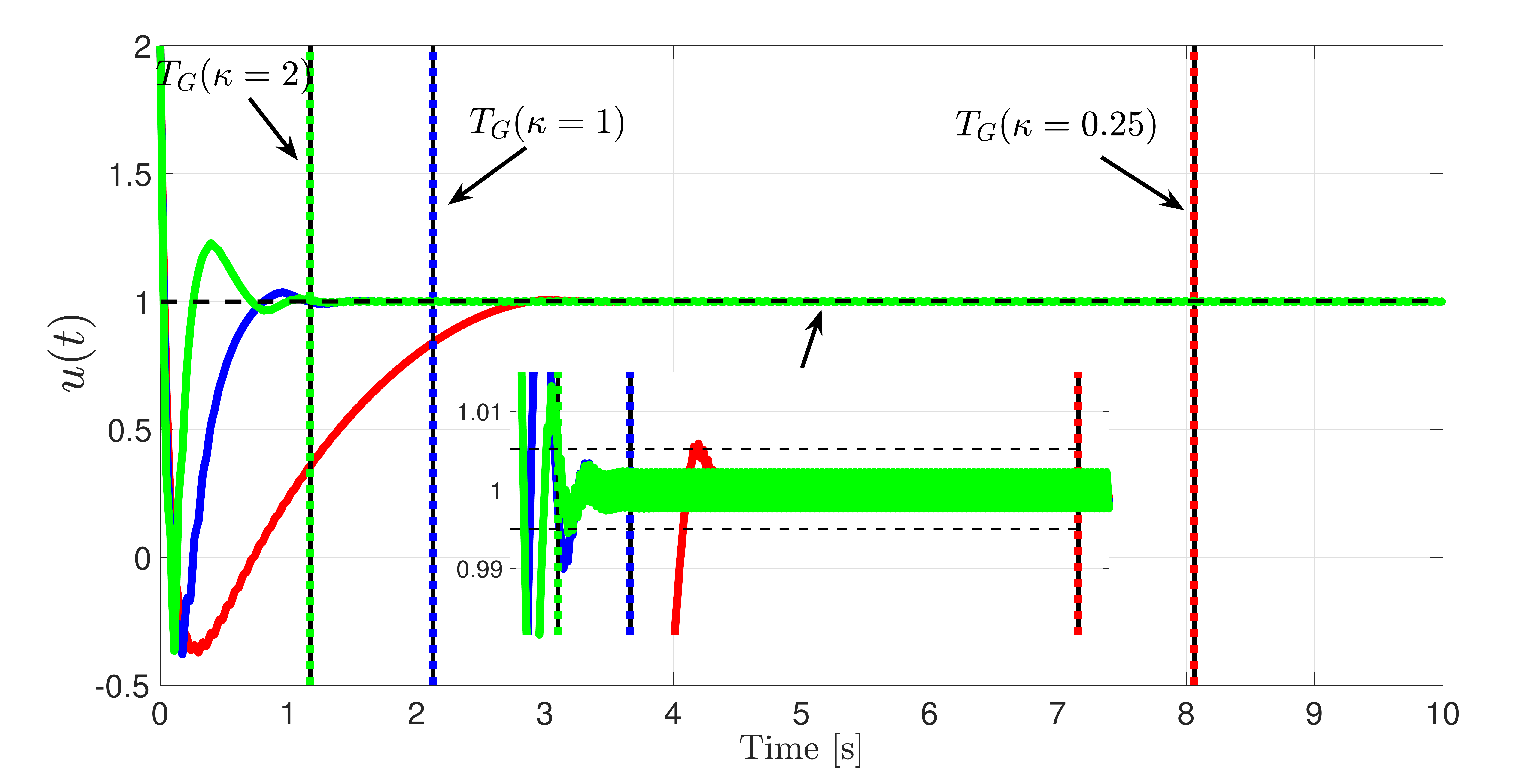}
    \caption{Evolution of $u$ for the cost function $\phi_{\kappa}$ of \eqref{cost_fixed} with three different parameters $\kappa\in\{0.25,1,2\}$ and with fixed parameters in the algorithm. Each trajectory corresponds to a different cost $\phi_{\kappa}$, and it converges to the set $[1-\nu,~1+\nu]$, with $\nu=0.003$, before the theoretical upper bound on the convergence time.}
    \label{figure1}
\end{figure}
From the proof of Theorem \ref{theorem1} we can observe that the fixed time convergence property, established in Step 1, is achieved by selecting admissible parameters $q_1$ and $q_2$ that also guarantee continuity of the vector field\footnote{The conference submission \cite{ACCpovedaKrstic} analyzed the FTGES dynamics using tools for discontinuous systems, e.g., Krasovskii regularizations. Since the Krasovskii regularization of a time-invariant continuous vector field is just the same vector field, said approach is still valid for the analysis. However, in this paper we are able to present a much simpler proof of the stability analysis of the FTGES by exploiting the continuity of the vector field \eqref{ES_dynamics1}.}, and by an appropriate orderly tuning of the parameters $(\varepsilon_2,a,\varepsilon_1)$. Moreover, since only positivity is required for the gain $k$, any fixed time $T_G^*>0$ can be assigned a priori by using equation \eqref{design_k}. Note, however, that $T_G^*$ depends on the  parameter $\kappa$ that characterizes the cost function in Assumption \ref{assumption_1}, which is assumed to be unknown. Therefore, in order to induce a particular fixed-time $T^*_G$, a conservative estimate of $\kappa$ should be used in order to tune the parameters of the FTGES.

\begin{remark}
The stability analysis of the FTGES relies on the convergence properties of the gradient flow \eqref{ES_dynamics_average_nominal_reduced}, which we use to qualitatively 
approximate the asymptotic behavior of the state $u$ in the FTGES. This is achieved by designing the FTGES dynamics \eqref{ES_dynamics1} to be in standard form for the application of singular perturbation theory \cite{Wang:12_Automatica,zero_order_poveda_Lina}, and by restricting a priori the state $\xi$ of the filter to the compact set $\eta\mathbb{B}$. Since, as shown in Step 2, the constant $\eta>0$ can be taken arbitrarily large to encompass any complete solution of practical interest, the restriction of $\xi$ to $\eta\mathbb{B}$ does not impose any practical constraint in the algorithm, and it is used only for the purpose of analysis. \QEDB
\end{remark}
\begin{remark}
Unlike the model-based fixed-time gradient dynamics of \cite{fixed_time} and \cite{Garg_Inequalities}, the FTGES dynamics are \emph{model-free} and only need measurements of the cost function $\phi$. Because of this, their stability properties are highly dependent on an appropriate \emph{ordered} tuning of the parameters $(\varepsilon_2,a,\varepsilon_1)$. Since standard averaging theory requires Lipschitz continuity of the vector field, the analysis of the FTGES dynamics cannot be carried out using standard averaging tools for smooth extremum seeking controllers such as in \cite{KrsticBookESC} or \cite{durrEbenbauer}. Instead, in this paper we used generalized averaging tools for non-smooth and hybrid systems \cite{Wang:12_Automatica,zero_order_poveda_Lina}.    \QEDB
\end{remark}

\begin{figure}[t!]
  \centering
    \includegraphics[width=0.49\textwidth]{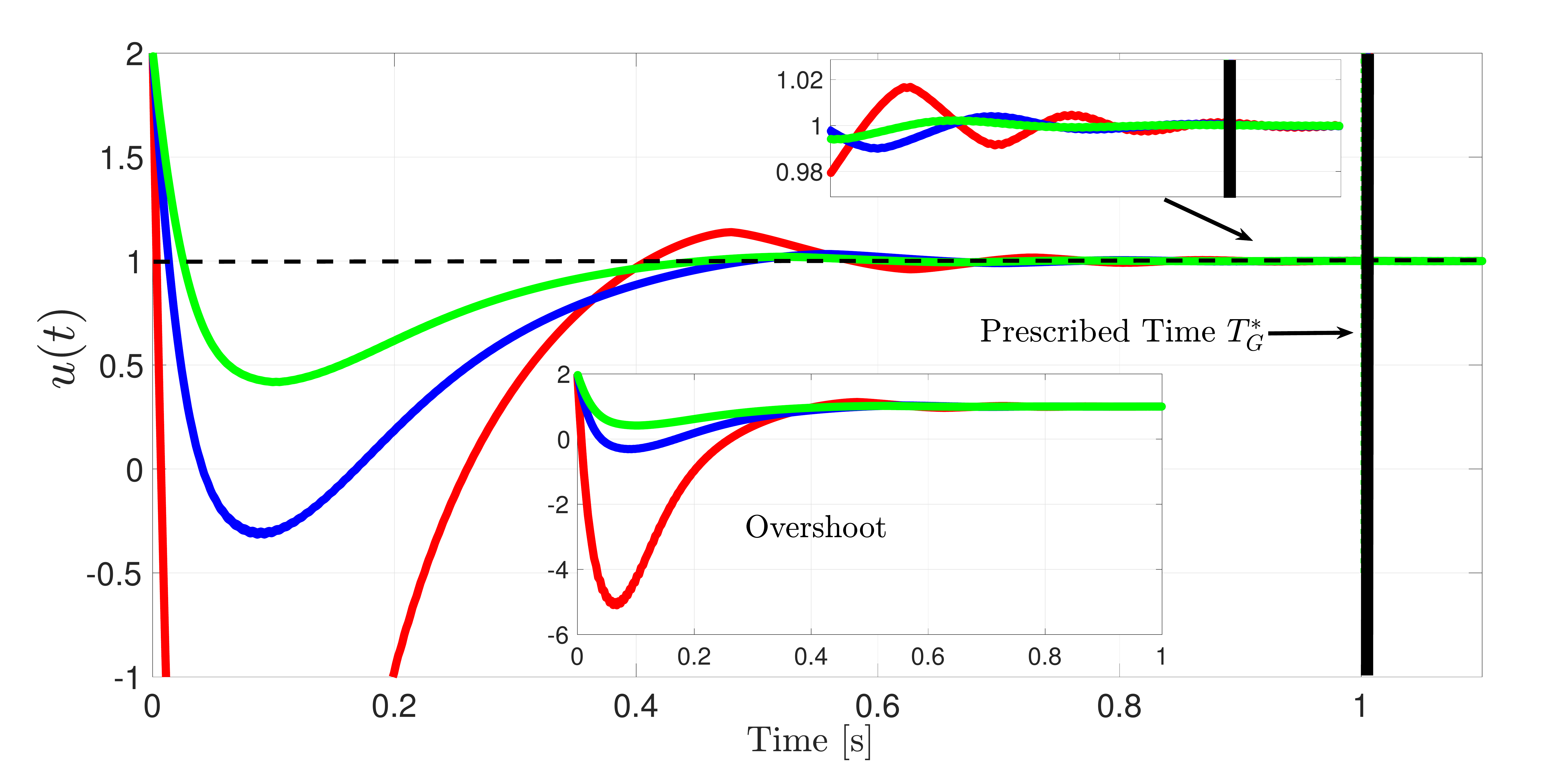}
    \caption{Evolution of $u$ for the cost function $\phi_{\kappa}$ of \eqref{cost_fixed} with three different values of $\kappa\in\{0.25,1,2\}$ and with fixed prescribed time $T_G^*=1.1$. Red line corresponds to $\kappa=0.25$, blue line corresponds to $\kappa=1$, and green line corresponds to $\kappa=2$. In this case, the parameters of the ES dynamics were tuned to guarantee convergence to the set $[1-\nu,~1+\nu]$, with $\nu=0.003$, before $T_G^*$ seconds. The insets show the overshoot and the convergence to the neighborhood of the optimizer.}
        \label{figure2}
\end{figure}

In order to illustrate the performance of the FTGES dynamics, we consider first the scalar cost function
\begin{equation}\label{cost_fixed}
\phi_{\kappa}(z)=\frac{\kappa}{2}(z-z^*)^2,~~~\kappa>0,
\end{equation}
which satisfies $\nabla^2\phi(z)=\kappa$. For the purpose of simulation, we initially consider the case where $q_1=3$,  $q_2=1.5$, $a=0.1$, $\varepsilon_1=0.02$, $\varepsilon_2=0.1$. We consider constants $\kappa\in\{0.25,1,2\}$ which generate the theoretical upper bounds $T_{0.25}=8.25,~T_{1}=2.06,~T_2=1.03$. Figure \ref{figure1} shows the behavior obtained under the FTGES dynamics. For each different $\kappa$ it can be seen that the trajectories converge to a small $\nu$-neighborhood of the optimal point in a finite time bounded by $T_G^*$. The controllers were tuned to induce this behavior with $\nu=0.003$.
On the other hand, we also consider the case when the finite time convergence is fixed at $T_G^*=1$. In this case, for all cost functions $\phi_{\kappa}$ the FTGES dynamics are tuned to guarantee that the convergence time is upper bounded by $1$. In order to achieve this property, the parameters were selected as $a=0.1$, $\varepsilon_1=0.001$, $\varepsilon_2=0.05$, $q_1=3$, $q_2=1.5$, and the gain $k$ was obtained via equation \eqref{design_k}. Figure \ref{figure2} shows the trajectories of $u$ for $\kappa\in\{0.25,1,2\}$. As expected all the trajectories converge to the set $[1-\nu,~1+\nu]$ in finite time bounded by $T_G^*=1$, where $\nu=0.005$. All the auxiliary states $\xi$ were initialized at zero.

\begin{figure}[t!]
  \centering
    \includegraphics[width=0.485\textwidth]{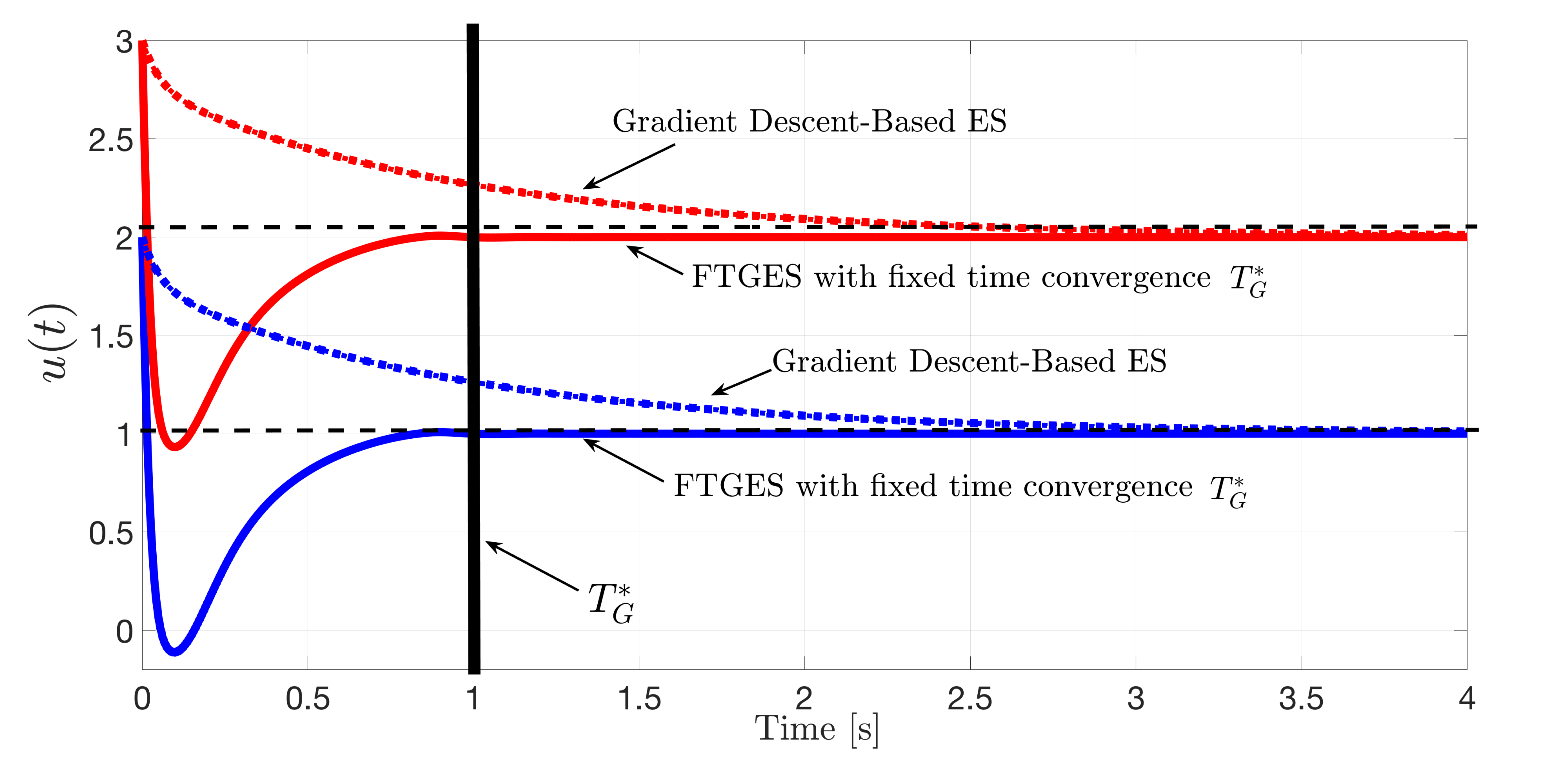}
    \caption{Trajectories of $u$ generated by the FTES dynamics applied to a multi-variable optimization problem. The dotted lines show the trajectories generated by the classic vanilla gradient descent based-ES.}
        \label{figure3}
\end{figure}

We finish this section by considering the multi-variable case, where the cost function has the form $\phi(z)=\frac{1}{2}z^\top Qz+b^\top z$.  We let $T_G^*=1$, and we tuned the parameters of the FTGES to generate trajectories that converge to a neighborhood of the optimal set in a finite time bounded by $T_G^*$. The matrix $Q\in\mathbb{R}^{2\times2}$ is selected as a diagonal positive definite matrix, and the vector $b$ is selected such that the unique minimizer of $\phi$ is $z^*=[1,2]^\top$. Figure \ref{figure3} shows the evolution of the solutions generated by the FTGES dynamics. The parameters were selected as $q_1=3$, $q_2=1.5$, $k=2.1$, $a=0.1$, $\varepsilon_1=0.001$, $\varepsilon_2=0.05$. It can be observed that the trajectories converge to a small neighborhood of the optimal point $z^*$ before $T_G^*$. In order to make a clear comparison with the standard ES dynamics, we have also have plotted the solutions of the classic gradient-based ES dynamics considered in \cite{Krstic2000,TanAndNesic2006Local}. 
%
%
%
\section{Newton-Based Fixed-Time Extremum Seeking}
\label{Sec_Newton}
While the FTGES dynamics \eqref{ES_dynamics1} are able to solve problem \eqref{main_problem} in a fixed time $T_G^*$, as shown in equations \eqref{gamma1} and \eqref{fixed_time1} the value of $T_G^*$ depends on the unknown coefficient $\kappa$ that characterizes the cost function in Assumption \ref{assumption_1}. This dependence was further illustrated in Figure \ref{figure1}, where the FTGES dynamics with identical parameters were applied to three cost functions with different coefficients $\kappa$. Since the standing assumption in extremum seeking is that $\phi$ is unknown, absence of knowledge of the coefficient $\kappa$ presents a challenge for the tuning of the FTGES.  Motivated by this limitation, in this section we now present a Fixed Time Newton-based Extremum Seeking  algorithm that removes the dependence on the coefficient $\kappa$
in the upper bound of the convergence time, thus facilitating the tuning of the control parameters.
\subsection{Qualitative Properties of the Cost Function}
For the Newton-based fixed time extremum seeking dynamics (FTNES) we consider the following assumption on the cost function $\phi$.
\begin{assumption}\label{assumption_2}
The function $\phi:\mathbb{R}^n\to\mathbb{R}$ is twice continuously differentiable, the Hessian matrix $\nabla^2 \phi(z)$ is positive definite for all $z\in\mathbb{R}^n$, and the norm $|\nabla \phi|$ has bounded level sets.  \QEDB
\end{assumption}
Since by assumption $\inf_{z\in\mathbb{R}^n}\phi(z)>-\infty$, the function $\phi$ attains a minimum. Thus, since $\nabla^2 \phi(z)>0$ implies strict convexity, the minimizer $z^*$ is unique and $\nabla \phi(z)=0$ if and only if $z=z^*$. Since there are strictly convex functions that do not satisfy the PL inequality \eqref{condition2}, the FTNES dynamics can be applied to cost functions $\phi$ that do not satisfy Assumption \ref{assumption_1}. 
\subsection{Newton-based Fixed-Time Dynamics}
In order to solve problem \eqref{main_problem} for functions satisfying Assumption \ref{assumption_2}, we now consider an extremum seeking controller with state $(u,\xi_1,\xi_2,\mu)\in\mathbb{R}^n\times\mathbb{R}^n\times\mathbb{R}^{n\times n}\times\mathbb{R}^{2n}$, evolving on the set 
\begin{equation}\label{flowset1N}
(u,\xi_1,\xi_2,\mu)\in C:=\mathbb{R}^n\times \eta\mathbb{B}\times \eta\mathbb{B} \times \mathbb{S}^n,
\end{equation}
with dynamics
\begin{align}\label{ES_dynamics1N}
\left(\begin{array}{c}
\dot{u}\\\\
\dot{\xi}_1\\\\
\dot{\xi}_2\\\\
\dot{\mu}
\end{array}\right)=-\left(\begin{array}{c}
k\xi_1\left(\dfrac{\xi_2}{|\xi_2|^{\alpha_1}}+\dfrac{\xi_2}{|\xi_2|^{\alpha_2}}\right)\\
\dfrac{1}{\varepsilon_2}\Big(\xi_1F_H(\phi,\mu)\xi_1-\xi_1\Big),\\
\dfrac{1}{\varepsilon_2}\Big(\xi_2-F_G(\phi,\mu)\Big)\\
\dfrac{2\pi}{\varepsilon_1}\mathcal{R}_{\kappa}\mu,
\end{array}\right).
\end{align}
where the right-hand side of $\dot{u}$ is defined as zero whenever $\xi_2=0$.  The parameters $(\alpha_1,\alpha_2)$ are again given by \eqref{alphaconstants}, the mapping $F_G$ is given by \eqref{mappingG}, the input $z$ and the function $M$ are given by \eqref{input}, and the dynamic oscillator is the same of \eqref{ES_dynamics1}. 

The FTNES has an extra state $\xi_1$ with dynamics depending on the mapping $F_H$, defined as
\begin{equation*}
F_H(\phi,\mu):=\phi(z)N(\mu),
\end{equation*}
where
\begin{equation}\label{excitation_signals}
N(\mu):=\left[\begin{array}{cccc}
N_{11} & N_{12} & \ldots & N_{1n}\\
N_{21} & N_{22} & \ldots & N_{2n}\\
\vdots  & \vdots & \vdots & \vdots\\
N_{n1} & N_{n2} & \ldots & N_{nn}
\end{array}\right],
\end{equation}
and where the entries $N_{ij}$ satisfy $N_{ij}=N_{ji}$, as well as the following conditions:
\begin{align*}
N_{ij}&=\frac{16}{a^2}\left(\tilde{\mu}_i^2-\frac{1}{2}\right),~~~~~\forall~i=j,\\
N_{ij}&=\frac{4}{a^2}\tilde{\mu}_i\tilde{\mu}_j,~~~~~~~~~~~~~~\forall~i\neq j.
\end{align*}
where $\tilde{\mu}=\mathcal{D}\mu$.

\begin{figure}[t!]
  \centering
    \includegraphics[width=0.4\textwidth]{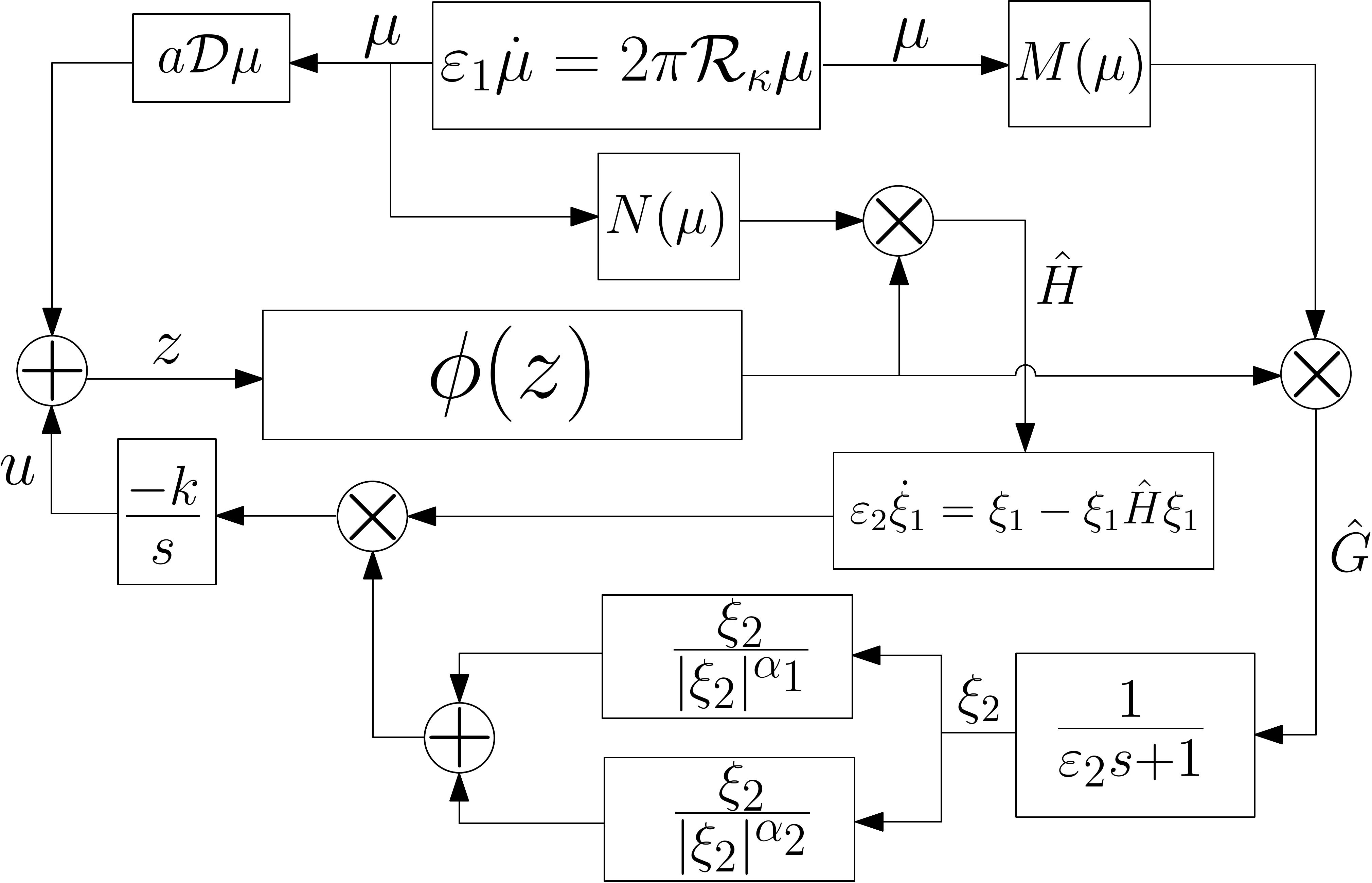}
    \caption{Scheme of the Fixed-Time Newton-based Extremum Seeking (FTNES) algorithm for a static map $\phi$.}
    \label{figureNewton}
\end{figure}

The structure of the NFTES is similar to the multi-variable Newton-based extremum seeking controller considered in \cite{Newton}, where, on average, the state $\xi_2$ serves as an estimation of the gradient, and the state $\xi_1$ approximates the inverse of the Hessian $\nabla^2\phi(z)$.  However, the NFTES dynamics have two main differences: a) for admissible parameters $(q_1,q_2)$, the dynamics of $u$ are continuous but not Lipschitz continuous, and they aim to approximate a Newton-based flow with fixed-time convergence properties \cite{fixed_time,Garg_Inequalities} instead of the standard Newton-based flow $\dot{x}=-\nabla^2\phi(x)^{-1}\nabla\phi(x)$ considered in \cite{Newton}; b) the excitation signals $N(\cdot)$ and $M(\cdot)$ are both generated by a time-invariant linear oscillator, which facilitates the analysis of the algorithm via averaging theory for non-smooth time-invariant dynamical systems \cite{Wang:12_Automatica}. A scheme illustrating the FTNES dynamics is shown in Figure \ref{figureNewton}.
\begin{remark}
In equation \eqref{ES_dynamics1N} the state $\xi_1$ corresponds to an $n\times n$ matrix. Therefore, the dynamics of $\xi_1$ must be understood as a matrix differential equation. This notation is used to simplify our presentation, and it is consistent with the standard Newton-based extremum seeking controllers of \cite{Newton}. There is no loss of generality in analyzing these dynamics in vectorial form, such as in \cite{NewtonESC}. \QEDB
\end{remark}
Since admissible parameters $(q_1,q_2)$ guarantee that $1-\alpha_1>0$ and $1-\alpha_2>0$, the dynamics of $u$ in \eqref{ES_dynamics1} also satisfy
\begin{align*}
\left|\lim_{\xi_2\to0}\frac{\xi_1\xi_2}{|\xi_2|^{\alpha_1}}\right|=\lim_{\xi_2\to0}\left|\frac{\xi_1\xi_2}{|\xi_2|^{\alpha_1}}\right|&\leq \lim_{\xi_2\to0} \left| \xi_1\right| \frac{|\xi_2|}{|\xi_2|^{\alpha_1}}\\
&\leq \left| \xi_1\right|\lim_{\xi_2\to0} |\xi_2|^{1-\alpha_1}\\
&=0.
\end{align*}
and
\begin{align*}
\left|\lim_{\xi_2\to0}\frac{\xi_1\xi_2}{|\xi_2|^{\alpha_2}}\right|=\lim_{\xi_2\to0}\left|\frac{\xi_1\xi_2}{|\xi_2|^{\alpha_2}}\right|&\leq \lim_{\xi_2\to0} \left| \xi_1\right| \frac{|\xi_2|}{|\xi_2|^{\alpha_2}}\\
&\leq \left| \xi_1\right|\lim_{\xi_2\to0} |\xi_2|^{1-\alpha_2}\\
&=0.
\end{align*}
Therefore, admissible parameters guarantee continuity of the FTNES dynamics. Moreover, since the algorithm has the same tunable parameters $(q_1,q_2,k,\varepsilon_1,\varepsilon_2,a)$ of the FTGES algorithm \eqref{ES_dynamics1}, under an appropriate tuning of the parameters a multi-time scale behavior can be induced in the closed-loop system. 

\subsection{Main Result}
In order to state the main convergence result for the FTNES dynamics, for each pair of admissible parameters $(q_1,q_2)$ we now define the positive constants
\begin{align*}
\tilde{\gamma}_1:=2^{\frac{\alpha_1}{2}},~~~~\tilde{\gamma}_2:=2^{\frac{\alpha_2}{2}},
\end{align*}
and the upper bound
\begin{equation}\label{upper_bound}
T_N^*:=\frac{1}{k}\left(\frac{\tilde{\gamma}_1}{\alpha_1}-\frac{\tilde{\gamma}_2}{\alpha_2}\right),
\end{equation}
where $(\alpha_1,\alpha_2)$ are defined as in \eqref{alphaconstants} and $k>0$. For admissible parameters $(q_1,q_2)$, the term inside the parentheses of \eqref{upper_bound} is positive. Moreover, $T_N^*$ does not depend on any parameter of the cost function $\phi$. Thus, for each desired $T_N^*>0$ one can always satisfy equation \eqref{upper_bound} by choosing admissible parameters $(q_1,q_2,k)$ with 
\begin{equation}\label{choiceC}
k=\frac{1}{T_N^*}\left(\frac{\tilde{\gamma}_1}{\alpha_1}-\frac{\tilde{\gamma}_2}{\alpha_2}\right).
\end{equation}
The following Theorem states the main convergence result for the FTNES dynamics:
\begin{thm}\label{main_theorem2}
Consider the NFTES dynamics and suppose that Assumptions \ref{assumptionFTGES} and \ref{assumption_2} hold. Then, for admissible parameters $(q_1,q_2)$, $k>0$ and each $\nu>0$, there exists $\eta>0$ and $\varepsilon_2^*>0$ such that for each $\varepsilon_2\in(0,\varepsilon_2^*)$ there exists $a^*>0$ such that for each $a\in(0,a^*)$ there exists $\varepsilon_1^*>0$ such that for each $\varepsilon_1\in(0,\varepsilon_1^*)$ there exists a neighborhood $\mathcal{N}$ of $p^*:=(z^*, \nabla^2 \phi(z^*)^{-1},\nabla\phi(z^*))$ such that every solution with $(u(0),\xi(0),\mu(0))\in \mathcal{N}\times\mathbb{S}^n$ is defined for all time $t\geq0$, and satisfies $|z(t)-z^*|\leq \nu$ for all $t\geq T_N^*$. \QEDB
\end{thm}
\textsl{Proof:} The proof of Theorem \ref{main_theorem2} is similar to the proof of Theorem \ref{theorem1}. We start with the next auxiliary Lemma which follows by direct computation.
\begin{lemma}\label{integrals2}
Under item 3) of Assumption \ref{assumptionFTGES} there exists a $T>0$ such that for all solutions $\mu$ of the linear oscillator \eqref{oscillator} with $\varepsilon_1=1$ the following holds:
\begin{align*}
&\frac{1}{\ell T}\int_{0}^{\ell T}\tilde{\mu}_i(s)ds=0,~~\forall~~i\in\{1,2,\ldots,n\}.\\
&\frac{1}{\ell T}\int_{0}^{\ell T}\tilde{\mu}_i(s)^2ds=\frac{1}{2},\\
&\frac{1}{\ell T}\int_{0}^{\ell T}\tilde{\mu}_{i}(s)\tilde{\mu}_j(s)ds=0,~~\forall~~i\neq j,\\
&\frac{1}{\ell T}\int_{0}^{\ell T}\tilde{\mu}_{i}(s)^2N_{ii}(s)ds=\frac{1}{8},~~\forall~~i\in\{1,2,\ldots,n\},\\
&\frac{1}{\ell T}\int_{0}^{\ell T}\tilde{\mu}_{i}(s)^2N_{jj}(s)ds=0,~~\forall~~i\neq j,\\
&\frac{1}{\ell T}\int_{0}^{\ell T}\tilde{\mu}_{i}(s)^2N_{ij}(s)ds=0,~~\forall~~i\neq j,\\
&\frac{1}{\ell T}\int_{0}^{\ell T}\tilde{\mu}_{i}(s)N_{ij}(s)ds=0,~~~\forall~~i\neq j,\\
&\frac{1}{\ell T}\int_{0}^{\ell T}\tilde{\mu}_{i}(s)N_{ii}(s)ds=0,~~\forall~~i\in\{1,2,\ldots,n\},\\
&\frac{1}{\ell T}\int_{0}^{\ell T}\tilde{\mu}_{i}(s)N_{jj}(s)ds=0,~~\forall~~i\neq j,\\
&\frac{1}{\ell T}\int_{0}^{\ell T}\tilde{\mu}_{i}(s)\tilde{\mu}_{j}(s)N_{ij}(s)ds=\frac{1}{4},~~\forall~~i\neq j,\\
&\frac{1}{\ell T}\int_{0}^{\ell T}\tilde{\mu}_{i}(s)\tilde{\mu}_{j}(s)N_{ii}(s)ds=0,~~\forall~~i\neq j.
\end{align*}
for all $\ell\in\mathbb{Z}_{\geq 1}$, where $\tilde{\mu}=\mathcal{D}\mu$.
  \QEDB
\end{lemma}
%
%
%
%
Lemma \ref{integrals2} is instrumental for the application of averaging theory in order to obtain estimations of $\nabla\phi$ and $\nabla^2 \phi$ of order $O(a)$. In particular, by taking again a Taylor expansion of $\phi(u+a\tilde{\mu})$ around $u$ for small values of $a$, and retaining the second order terms, we obtain:
\begin{equation*}
\phi(u+a\tilde{\mu})=\phi(u)+a\tilde{\mu}^\top \nabla \phi(u)+\frac{a^2}{2}\tilde{\mu}^\top\nabla^2 \phi(u)\tilde{\mu}+O(a^3).
\end{equation*}
Using this expansion, and the definitions of the mappings $M$, $N$, $F_G$, and $F_H$, as well as the integrals of Lemma \ref{integrals2}, we obtain the following average functions, where the average is taken with respect to the solutions of the oscillator, i.e., keeping $(u,\xi_1,\xi_2)$ constant:
\begin{equation*}
\frac{1}{\ell T}\int_{0}^{\ell T}F_{G}(\phi(s),\mu(s))ds=\nabla\phi(u)+O(a),
\end{equation*}
and
\begin{equation*}
\frac{1}{\ell T}\int_{0}^{\ell T}F_{H}(\phi(s),\mu(s))ds=\nabla^2\phi(u)+O(a).
\end{equation*}
Therefore, the average dynamics of \eqref{ES_dynamics1N} are given by
\begin{align}\label{NES_dynamics1}
\left(\begin{array}{c}
\dot{u}^a\\\\
\dot{\xi}^a_1\\\\
\dot{\xi}^a_2\\
\end{array}\right)=-\left(\begin{array}{c}
k\xi^a_1\left(\dfrac{\xi^a_2}{|\xi^a_2|^{\alpha_1}}+\dfrac{\xi^a_2}{|\xi^a_2|^{\alpha_2}}\right)\\
\dfrac{1}{\varepsilon_2}\Big(\xi^a_1\nabla^2\phi(u^a)\xi^a_1-\xi^a_1+O(a\xi_1^{2a})\Big),\\
\dfrac{1}{\varepsilon_2}\Big(\xi^a_2-\nabla \phi(u^a)+O(a)\Big)
\end{array}\right).
\end{align}
System \eqref{NES_dynamics1} is a perturbed version of the nominal average system
\begin{align}\label{NES_dynamics2}
\left(\begin{array}{c}
\dot{u}^a\\\\
\dot{\xi}^a_1\\\\
\dot{\xi}^a_2\\
\end{array}\right)=-\left(\begin{array}{c}
k\xi^a_1\left(\dfrac{\xi^a_2}{|\xi^a_2|^{\alpha_1}}+\dfrac{\xi^a_2}{|\xi^a_2|^{\alpha_2}}\right)\\
\dfrac{1}{\varepsilon_2}\Big(\xi^a_1\nabla^2\phi(u^a)\xi^a_1-\xi^a_1\Big),\\
\dfrac{1}{\varepsilon_2}\Big(\xi^a_2-\nabla \phi(u^a))\Big)
\end{array}\right).
\end{align}
Thus, we analyze again the stability properties of system \eqref{NES_dynamics1} by studying the properties of system \eqref{NES_dynamics2}, and using standard robustness results for perturbed ODEs with a continuous right-hand side, e.g., \cite[Thm. 7.21]{HDS}.

When $\varepsilon_2$ is sufficiently small, system \eqref{NES_dynamics2} is a singularly perturbed system. The boundary layer dynamics are obtained in the $\tau=t/\varepsilon_2$ time scale by setting $\varepsilon_2=0$, which leads to $\frac{\partial u^a}{\partial \tau}=0$ and
\begin{subequations}\label{ES_dynamics3}
\begin{align}
\frac{\partial \xi^a_1}{\partial \tau}&=\xi^a_1-\xi^a_1\nabla^2\phi(u^a)\xi^a_1,\label{hessian_dynamics3}\\
\frac{\partial \xi^a_2}{\partial \tau}&=-\xi^a_2 +\nabla \phi(u^a).\label{gradient_dynamics3}
\end{align}
\end{subequations}
For fixed values of $u^a$, the stability properties of system \eqref{ES_dynamics3} can be analyzed as follows: denote by $r_u:=\nabla \phi(u^a)$ and $H_u:=\nabla^2\phi(u^a)$, and consider the errors $\tilde{\xi}^a_1=\xi^a_1-H_u^{-1}$ and $\tilde{\xi}^a_2=\xi^a_2-r_u$. The error boundary layer dynamics are then given by the following decoupled equations:
\begin{subequations}
\begin{align}
\dot{\tilde{\xi}}^a_1&=-\tilde{\xi}^a_1H_u\left(\tilde{\xi}^a_1+H_u^{-1}\right)\label{filter1}\\
\dot{\tilde{\xi}}^a_2&=-\tilde{\xi}^a_2.\label{filter2}
\end{align}
\end{subequations}
System \eqref{filter2} has the origin globally exponentially stable, and system \eqref{filter1} has the origin locally exponentially stable given that its linearization around the origin has the Jacobian $-I$, see \cite[pp. 1761]{Newton}. Therefore, for each fixed $u^a$, the boundary layer dynamics render the quasi-steady state mapping
\begin{equation*}
 \xi^*=\left(\nabla^2 \phi(u^a)^{-1},\nabla \phi(u^a)\right),
\end{equation*}
locally exponentially stable, uniformly on $u^a$. Using this stability property we can now obtained the reduced dynamics associated to system \eqref{NES_dynamics2} by substituting $\xi^*$ in the dynamics of $u^a$:
\begin{align}\label{slow_dynamics}
\dot{u}_r=-k\nabla^2\phi(u_r)^{-1}\left(\frac{\nabla\phi(u_r)}{|\nabla\phi(u_r)|^{\alpha_1}}+\frac{\nabla\phi(u_r)}{|\nabla\phi(u_r)|^{\alpha_2}}\right).
\end{align}
Following the ideas of \cite{fixed_time}, we can analyze the stability properties of system \eqref{slow_dynamics} by considering the Lyapunov function
\begin{equation*}
V_H(u_r)=\frac{1}{2}|\nabla \phi(u_r)|^2,
\end{equation*}
which, under Assumption \ref{assumption_2}, is radially unbounded and positive definite with respect to the point $\{z^*\}$. The derivative of $V_H$ along the solutions of \eqref{slow_dynamics} satisfies
\begin{equation}\label{Lyapunov}
\dot{V}_H(u_r)=-k\rho_1V_H(u_r)^{\chi_1}-k\rho_2V_H(u_r)^{\chi_2}<0,
\end{equation}
for all $u_r\neq z^*$, where
\begin{align*}
\rho_1=2^{\chi_1}>0,~~~\rho_2=2^{\chi_2}>0
\end{align*}
and 
\begin{align*}
\chi_1=\frac{2-\alpha_1}{2}\in (0.5,1),~~~\chi_2=\frac{2-\alpha_2}{2}>1.
\end{align*}
Therefore, system \eqref{slow_dynamics} renders the point $z^*$ UGAS and by  \cite[Lemma 1]{Fixed_timeTAC}, the Lyapunov function evaluated along the solutions of \eqref{slow_dynamics} satisfies $V_H(u_r(t))=0$ for all $t\geq T_N^*$, where $T_N^*$ is given by \eqref{upper_bound}. Thus the convergence of $u_r$ to $z^*$ occurs in a finite time bounded by $T_N^*$. UGAS of  $\{z^*\}$ implies the existence of a class $\mathcal{K}\mathcal{L}$ function $\beta'_u$, constructed as \eqref{KLexplicit}, such that all solutions of \eqref{slow_dynamics} satisfy the bound
\begin{equation*}
|u_r(t)-z^*|\leq \beta'_u(|u_r(0)-z^*|,t),
\end{equation*}
for all $t\geq0$, and where $\beta'_u(r,t)=0$ for all $t\geq T^*$ and all $r>0$. From here, the stability analysis follows the exact same steps of the proof of Theorem \ref{theorem1}, using \cite[Thm.2]{Wang:12_Automatica} and \cite[Thm. 7.21]{HDS} successively to link the stability properties of the original dynamics \eqref{ES_dynamics1} and its average perturbed dynamics \eqref{NES_dynamics1}. Local existence of complete solutions from the neighborhood $\mathcal{N}$ follows directly by the local stability properties of \eqref{NES_dynamics2} and the closeness of solutions between \eqref{NES_dynamics2} and \eqref{ES_dynamics1N}.  \null\hfill\null $\blacksquare$

\subsection{Discussion and Numerical Examples}
\label{sec_NumNewton}
The convergence result of Theorem \ref{main_theorem2} is local with respect to the initial conditions, and practical with respect to the parameters $(\varepsilon_1,\varepsilon_2,a)$ and the neighborhood $\{z^*\}+\nu\mathbb{B}$. However, unlike existing results in the literature, an upper bound on the convergence time of the FTNES dynamics can be prescribed a priori by selecting admissible parameters $(q_1,q_2)$ and a gain $k$ that satisfies equation \eqref{choiceC}, and by initializing the states $(\xi_1,\xi_2)$ in a pre-defined compact set $\eta\mathbb{B}$. This represents a clear advantage in comparison to the existing Newton-based extremum seeking algorithms. 
\begin{remark}
The local convergence result of Theorem \ref{main_theorem2} is due to the existence of multiple equilibria in the dynamics
\begin{equation}\label{estimator_Hessian2}
\dot{\xi}_1=\Big(\xi_1-\xi_1\nabla^2\phi(z) \xi_1\Big),
\end{equation}
which emerge in the nominal average system \eqref{NES_dynamics2}. Similar local results emerge in Newton-based ESCs with asymptotic convergence properties; see for instance \cite{Newton}. While it is possible to design Newton-based ESCs with semi-global practical asymptotic stability results by computing the vector $\nabla^2\phi(x)^{-1}\nabla\phi(x)$ using the singular perturbation approach presented in \cite[Sec. 3]{NewtonSemiglobal}, see for instance \cite{NewtonESC}, said approach cannot be used in this case since it will generate learning dynamics for the state $u$ with discontinuous vector fields that are not locally bounded when $\nabla\phi=0$. \QEDB
\end{remark}
\begin{figure}[t!]
 \begin{centering}
  \includegraphics[width=0.48\textwidth]{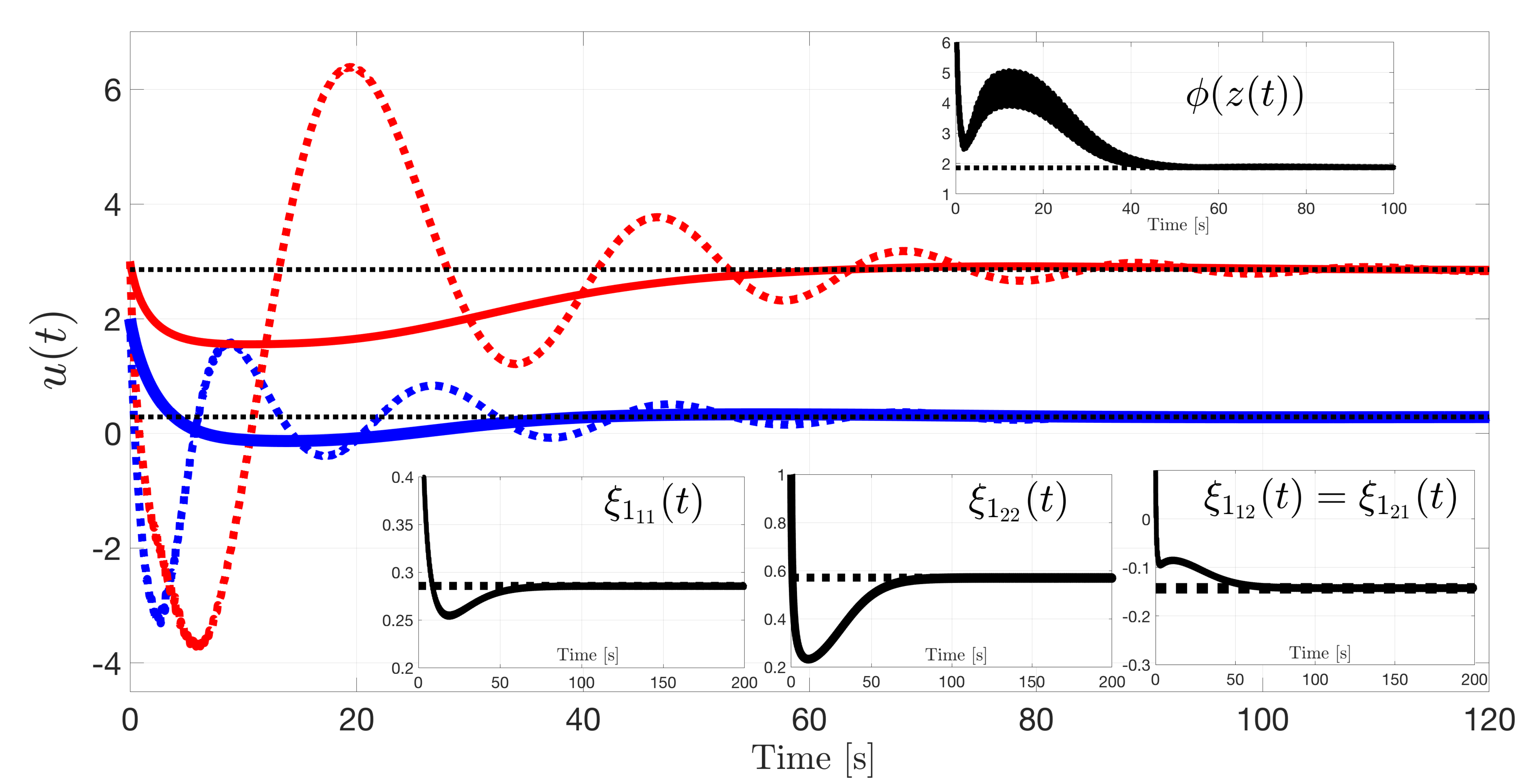}
  \caption{Evolution in time of the state $u$. Blue line corresponds to $u_1$, and red line corresponds to $u_2$. The dotted lines correspond to the trajectories generated by the traditional Newton-based ESC of \cite{Newton}. The solid lines correspond to the trajectories generated by the NFTES. The inset shows the trajectories of the state $\xi_1$ and the cost function $\phi(z)$.}
  \label{Fig11}
 \end{centering}
\end{figure}
To illustrate the performance of the NFTES, and to highlight the differences with respect to the the standard Newton-based extremum seeking controller of \cite{Newton}, consider the quadratic function
\begin{equation*}
\phi(z)=\frac{1}{2} z^\top H z+b^\top z+c,
\end{equation*}
which satisfies
\begin{align*}
\nabla \phi(z)&=Hz+b,\\
\nabla^2\phi(z)&=H.
\end{align*}
The parameters of the cost function are selected as
\begin{equation*}
H=\left[\begin{array}{cc}
4 & 1 \\ 1 & 2
\end{array}\right],~~~b=\left[\begin{array}{c}
-4\\
-6
\end{array}\right],~~c=11.
\end{equation*}
The inverse of the Hessian matrix is given by
\begin{equation*}
H^{-1}=\left[\begin{array}{cc}
0.2857 & -0.1429 \\ -0.1429 & 0.5714
\end{array}\right],
\end{equation*}
and the function $\phi(z)$ has a global minimizer at the point
\begin{equation*}
z^*=-H^{-1}b=\left[\frac{2}{7},~\frac{20}{7}\right]^\top.
\end{equation*}
In order to find $z^*$ in fixed time, we implement the NFTES dynamics with parameters $a=0.1$, $\varepsilon_1=0.1$ and $\varepsilon_2=10$. The constants $(k,q_1,q_2)$ were selected as $\eta=50,$ $k=0.025$, $q_1=3$, and $q_2=1.5$,  which generates an upper bound $T_{N}^*$ approximately equal to $T_N^*=123.4$. We have also simulated the Newton-based ESC of \cite{Newton}, which has learning dynamics of the form $\dot{x}=-\xi_1\xi_2$. To obtain a smooth approximation of $H^{-1}$, we have also implemented an additional low-pass filter that receives as input $\xi_1$ and $\xi_2$, and which generates filtered outputs $\xi^f_1$ and $\xi^f_{2}$ that serve as inputs to the learning dynamics. As shown in \cite{Newton}, the incorporation of these filters does not affect the stability analysis of the algorithm. Figure \ref{Fig11} shows the trajectories generated by the NFTES dynamics as well as the trajectories of the standard Newton-based ESC. It can be observed that the NFTES dynamics exhibit a much better transient performance in terms of less oscillations and faster convergence time to a neighborhood of $z^*$. The insets show the evolution of the components of the state $\xi_1$, which correspond to the entries of $H^{-1}$, as well as the evolution in time of the cost function $\phi(z)$. As it can be observed, the trajectories of the entries of $\xi_1$ converge to a neighborhood of the true values of the entries of $H^{-1}$.

\begin{figure}[t!]
 \begin{centering}
  \includegraphics[width=0.48\textwidth]{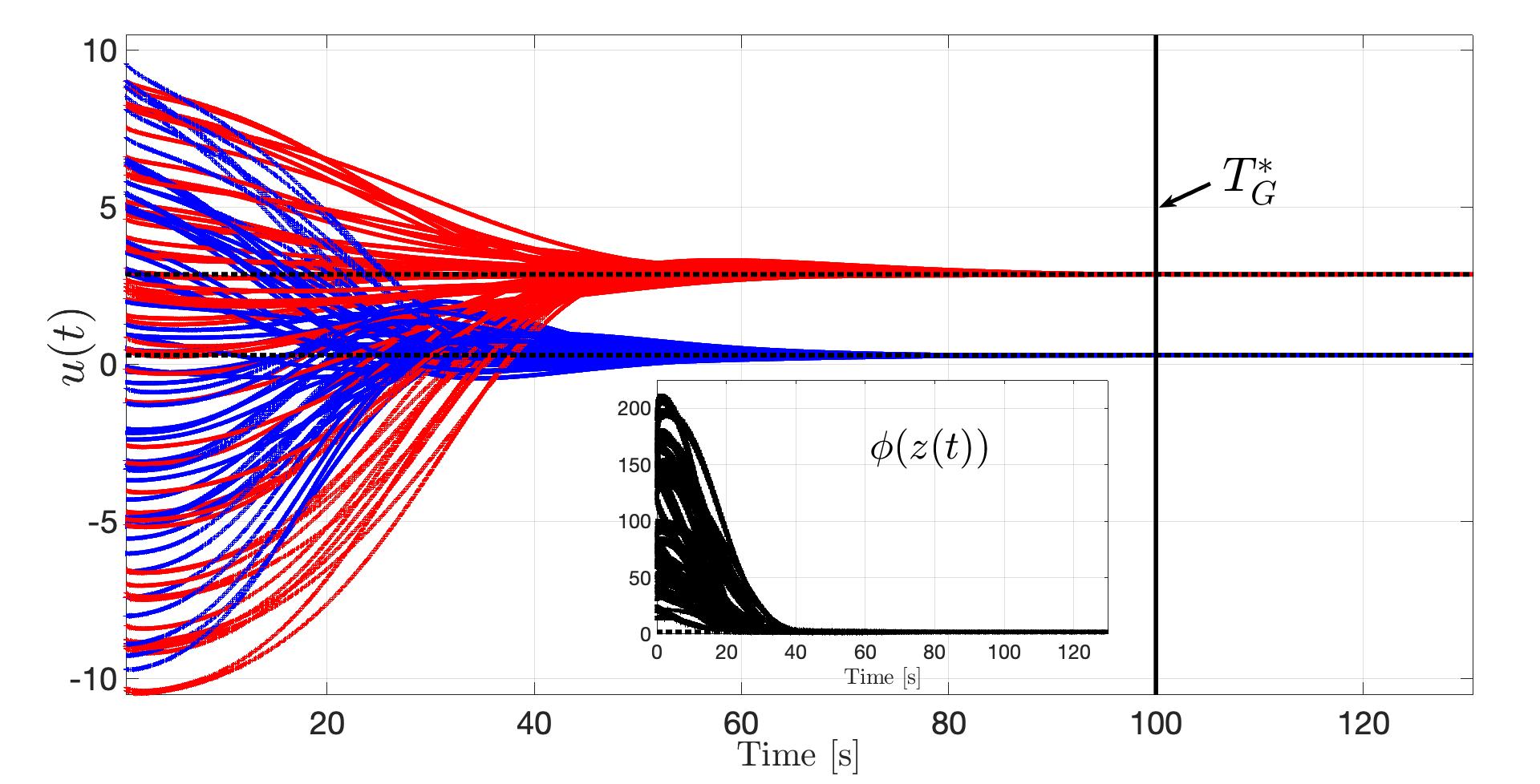}
  \caption{Evolution in time of several solutions $x$ initialized in a neighborhood of the optimizers. The solid black line indicates the upper bound $T^*$ on $T_c$.}
  \label{Fig6}
 \end{centering}
\end{figure}

To further illustrate the fixed-time convergence property of the NFTES dynamics, we have also simulated the case where the upper bound in the convergence time is fixed to $T_N^*=100$, which can be obtained in the NFTES dynamics by choosing $k=0.03085$, $q_1=3$, and $q_2=1.5$. Figure \ref{Fig6} shows the evolution in time of 50 different trajectories $u(t)$ initialized randomly in the set $[-10,10]\times[-10,10]$. The inset shows the evolution in time of the cost functions $\phi(z)$ along the trajectories of $z$. As it can be observed, the NFTES dynamics guarantee convergence to a small neighborhood of $z^*$ before the time $T^*$. 

%
%
%
%
%
It is important to note that in order to obtain the convergence result of Theorem \ref{main_theorem2}, the parameters $(a,\varepsilon_1,\varepsilon_2)$ need to be appropriately tuned. In particular, these parameters depend on the constants $(q_1,q_2,k)$, which, in turn, fix the bound $T_N^*$. Thus, smaller values of $T_N^*$ may require smaller values for $(a,\varepsilon_1,\varepsilon_2)$. Since $\varepsilon_2$ is related to the frequency of the oscillator, in order to implement in discrete time the NFTES algorithm for small values of $T_N^*$, one may need to use a sufficiently small step size to avoid aliasing issues. This indicates a potential tradeoff between achieving small finite convergence times and the computational complexity of the discretization of the algorithm. Also, larger sets $\eta\mathbb{B}$ may require smaller values of $(a,\varepsilon_1,\varepsilon_2)$. 

\section{$\varepsilon_0$-Fixed-Time Extremum Seeking for Dynamic Systems}
\label{Sec_Dynamic}
We now consider the more general extremum seeking problem, where the cost function $\phi$ corresponds to the steady state input-to-output mapping of a dynamical system. In particular, we consider the following dynamical system
\begin{subequations}\label{plant_dynamics}
\begin{align}
\dot{x}&=f(x,z)\label{plant_dynamics1}\\
y&=h(x,z)\label{output_plant},
\end{align}
\end{subequations}
where $x\in\mathbb{R}^p$ is the state of the system, $z\in\mathbb{R}^n$ is the input, $f$ is a continuous function characterizing the dynamics of the plant, and $h$ is an output function. We assume that the plant \eqref{plant_dynamics1} has already been stabilized such that the state $x$ evolves in a compact set $\Xi\subset\mathbb{R}^p$ for any input $z$ of interest. While this may look like a strong condition, it is a reasonable assumption given that we will consider plants \eqref{plant_dynamics} that have a quasi-steady state continuous manifold $z\mapsto \ell(z)$, and we will tune our controller to operate from particular compact sets $K_z$, which will generate uniformly bounded trajectories $z$ that will keep $x$ uniformly bounded.

In order to have a well-defined extremum seeking problem we also make the following stability assumption.
\begin{assumption}\label{assumption_plant}
There exists a continuous quasi-steady state manifold $\ell:\mathbb{R}^n\to\mathbb{R}^p$  such that for each compact set $K_z\subset\mathbb{R}^n$ the dynamics of the plant, with frozen input, given by
\begin{equation*}
(x,z)\in \Xi\times K_z,~\left\{\begin{array}{l}
\dot{x}=f(x,z)\\
\dot{z}=0
\end{array}\right.
\end{equation*}
render UGAS the compact set $M_{K_z}:=\{(x,z)\in\Xi\times K_z: x=\ell(z)\}$.\QEDB
\end{assumption}
The stability conditions of Assumption \ref{assumption_plant} are standard in extremum seeking problems, and they can be further relaxed to allow for set-valued quasi-steady state manifolds $\ell:\mathbb{R}^n\rightrightarrows\mathbb{R}^p$, see \cite{Poveda:16,black_box_CDC}. However, for simplicity we assume that $\ell$ is a continuous function, which allow us to avoid introducing extra definitions for set-valued maps.

Using the quasi-steady state manifold $\ell$ and the output \eqref{output_plant}, we define the quasi steady-state input-to-output mapping of the dynamical system \eqref{plant_dynamics} as
\begin{equation}\label{ssiom}
\phi(z):=h(\ell(z)).
\end{equation}
As before, we will assume that $\phi$ attains its minimum at some point $z^*\in\mathbb{R}^n$, i.e., $\inf_{z\in\mathbb{R}^n}\phi(z)>-\infty$.
\subsection{Closed-Loop Dynamics}
In order to optimize the steady-state input-to-output mapping of system \eqref{plant_dynamics} using the FTGES, we consider the following closed-loop system
\begin{align}\label{ES_dynamic}
\left(\begin{array}{c}
\dot{u}\vspace{0.1cm}\\
\dot{\xi}\vspace{0.1cm}\\
\dot{\mu}\vspace{0.1cm}\\
\dot{x}\vspace{0.1cm}
\end{array}\right)=\left(\begin{array}{c}
-k_1\left(\dfrac{\xi}{|\xi|^{\alpha_1}}+\dfrac{\xi}{|\xi|^{\alpha_2}}\right)\\
-k_2\big(\xi-F_G(y,\mu)\big)\\
-k_3\mathcal{R}_{\kappa}\mu\\
f(x,u+a\mathcal{D}u)
\end{array}\right),
\end{align}
evolving in the set
\begin{equation}\label{flowset11}
(u,\xi,\mu,x)\in C:=\mathbb{R}^n\times \eta\mathbb{B} \times \mathbb{S}^n\times \Xi,
\end{equation}
where $k_1:=\varepsilon_0 k$, $k_2:=\frac{\varepsilon_0}{\varepsilon_2}$, $k_3:=\frac{\varepsilon_02\pi}{\varepsilon_1}$, and where $\varepsilon_0$ is a new tunable parameter that satisfies $0<\varepsilon_0\ll 1$.
\begin{figure}[t!]
 \begin{centering}
  \includegraphics[width=0.48\textwidth]{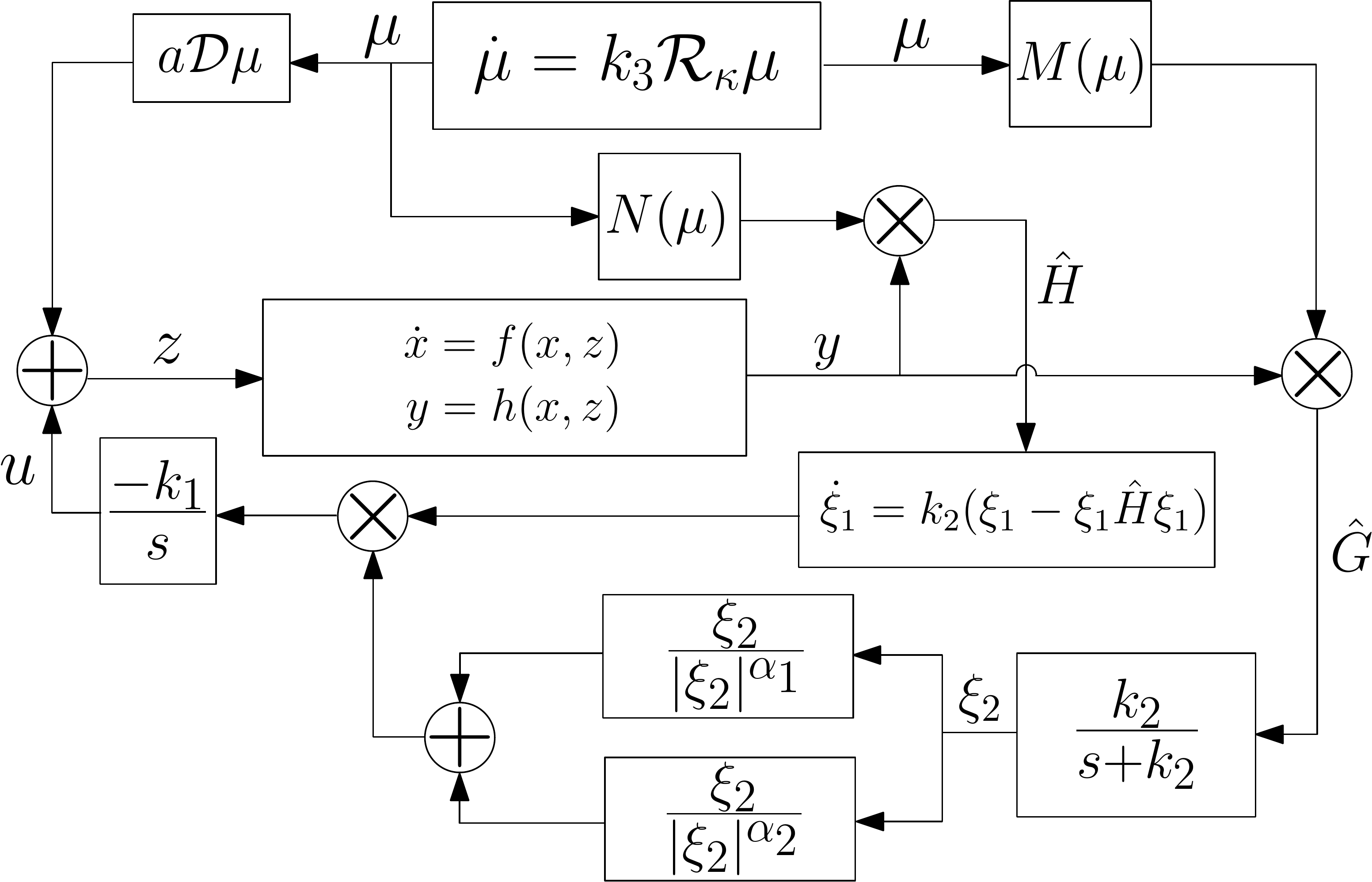}
  \caption{\small{Closed-loop system generated by the interconnection of the FTNES dynamics and the plant \eqref{plant_dynamics}.}}
  \label{Fig66}
 \end{centering}
\end{figure}
In contrast to the static case considered in Section \ref{Sec_Gradient}, the FTGES dynamics in \eqref{ES_dynamic} receive as input direct measurements of the output of the plant \eqref{plant_dynamics}. Since the rate of convergence of the state $x$ is unknown, it is not possible to prescribed a priori a convergence time for the complete closed-loop system operating in the original time scale $t$. Therefore, we study the convergence properties of system \eqref{ES_dynamic} in a new time scale $\tau:=t\varepsilon_0$. It then follows that $dt =d\tau/\varepsilon_0$, and the dynamics \eqref{ES_dynamic} in the $\tau$-time scale become
\begin{align}\label{ES_dynamicss}
\left(\begin{array}{c}
\dfrac{du}{d\tau}\vspace{0.1cm}\\
\dfrac{d\xi}{d\tau}\vspace{0.1cm}\\
\dfrac{d\mu}{d\tau}\vspace{0.1cm}\\
\dfrac{dx}{d\tau}
\end{array}\right)=\left(\begin{array}{c}
-k\left(\dfrac{\xi}{|\xi|^{\alpha_1}}+\dfrac{\xi}{|\xi|^{\alpha_2}}\right)\\
-\dfrac{1}{\varepsilon_2}\Big(\xi-F_G(y,\mu)\Big)\\
-\dfrac{2\pi}{\varepsilon_1}\mathcal{R}_{\kappa}\mu\\
\dfrac{1}{\varepsilon_0}f(x,u+a\mathcal{D}\mu)
\end{array}\right).
\end{align}
System \eqref{ES_dynamicss} is now in standard form for the application of singular perturbation theory. The boundary layer dynamics are obtained by setting $\varepsilon_0=0$ in \eqref{ES_dynamic}, which leads to $\dot{u}=0$, $\dot{\xi}=0$, $\dot{\mu}=0$, and the dynamics
\begin{equation}\label{bl_dynamic}
x\in\Xi,~~~\dot{x}=f(x,u+a\mathcal{D}\mu),
\end{equation}
which, by Assumption \ref{assumption_plant}, renders UGAS the quasi-steady state manifold $x=\ell(u+a\mathcal{D}\mu)$. Thus, the reduced system associated to the singularly perturbed system \eqref{ES_dynamicss} is given by
\begin{align}\label{ES_dynamicsss}
\left(\begin{array}{c}
\dfrac{du}{d\tau}\vspace{0.1cm}\\
\dfrac{d\xi}{d\tau}\vspace{0.1cm}\\
\dfrac{d\mu}{d\tau}\vspace{0.1cm}
\end{array}\right)=\left(\begin{array}{c}
-k\left(\dfrac{\xi}{|\xi|^{\alpha_1}}+\dfrac{\xi}{|\xi|^{\alpha_2}}\right)\\
-\dfrac{1}{\varepsilon_2}\Big(\xi-F_G(h(\ell(u+a\mathcal{D}\mu)),\mu)\Big)\\
-\dfrac{2\pi}{\varepsilon_1}\mathcal{R}_{\kappa}\mu
\end{array}\right).
\end{align}
Using the definition of $\phi$ in \eqref{ssiom} we can see that this system corresponds to the same system \eqref{ES_dynamics1}. From the proof of Theorem \ref{theorem1} we know that this system renders SGPAS as $(\varepsilon_2,a,\varepsilon_1)\to0^+$ the set $\mathcal{A}=\{z^*\}\times\eta\mathbb{B}\times\mathbb{S}^n$ with $\mathcal{K}\mathcal{L}$ bound $\beta_u$, where the parameters $(\varepsilon_2,a,\varepsilon_1)$ must be tuned orderly.  Since the boundary layer dynamics \eqref{bl_dynamic} are constrained to the set $\Xi$, the complete system \eqref{ES_dynamicss} renders SGPAS as $(\varepsilon_2,a,\varepsilon_1)\to0^+$ the set $\mathcal{A}=\{z^*\}\times\eta\mathbb{B}\times\mathbb{S}^n\times\Xi$ with $\mathcal{K}\mathcal{L}$ bound $\beta_u$, i.e., from compact sets of initial conditions and under suitable order tuning of the parameters of the controller, every solution of the closed-loop system generates trajectories $u$ that satisfy
\begin{equation*}
|u(\tau)-z^*|\leq \beta_u(|u(0)-z^*|,\tau)+\nu,
\end{equation*}
for all $t$ in the domain of the solutions, where $\nu$ can be made arbitrarily small, and where $\beta_u(|u(0)-z^*|,\tau)=0$ for all $\tau\geq T_G^*$, with $T_G^*$ given by \eqref{fixed_time1}. Therefore, in the $\tau$-time scale, the closed-loop system achieves fixed-time extremum seeking.  
\begin{figure*}[t!]
 \begin{centering}
  \includegraphics[width=0.75\textwidth]{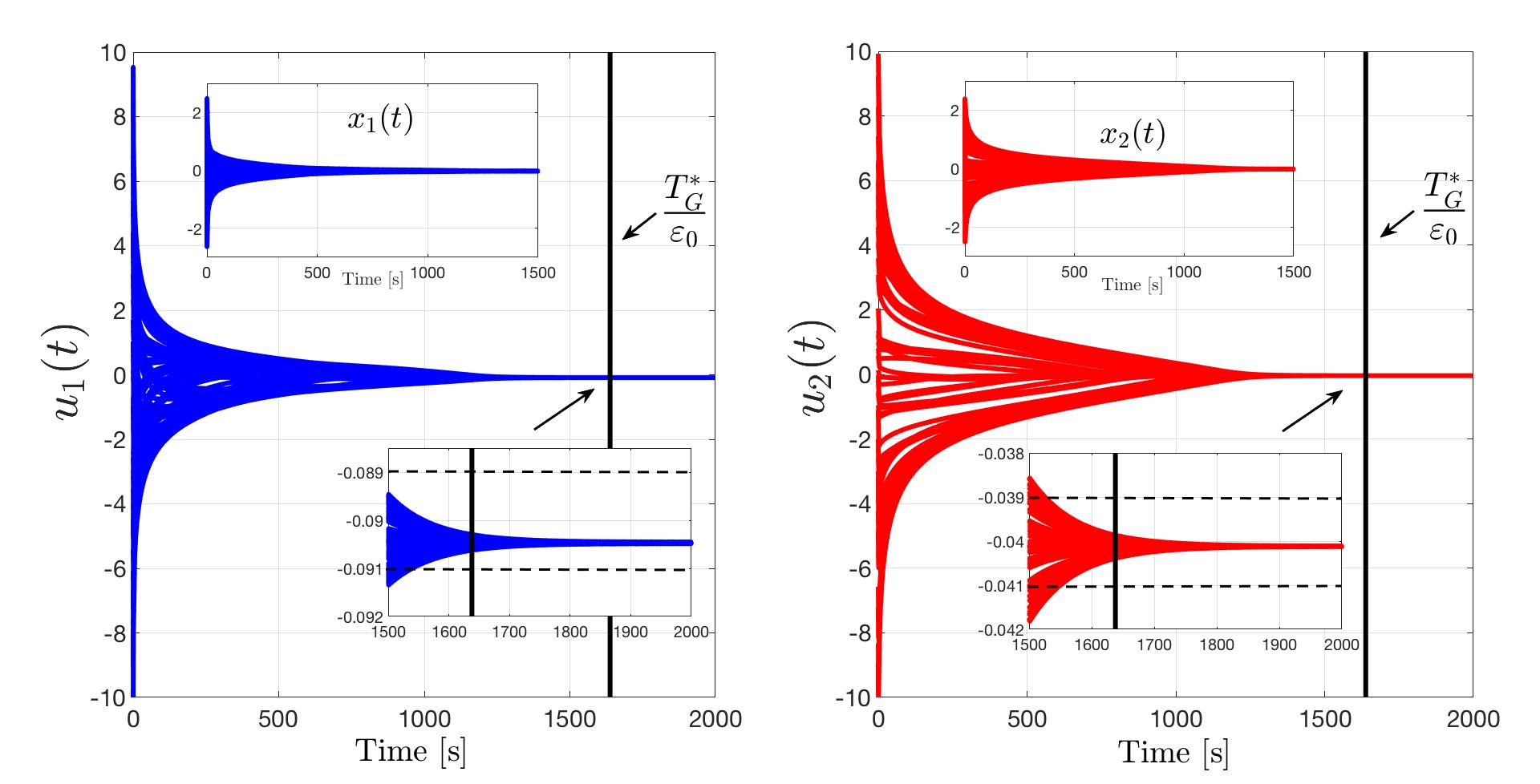}
  \caption{\small Time history of 70 different trajectories of the the states $u_1$ and $u_2$ generated by the FTGES from $70$ different initial conditions, and applied to the dynamic plant \eqref{dynamic_plant_example}. The trajectories were randomly initialized in the interval $[-10,10]$. The inset shows the trajectories of the states $x_1$ and $x_2$ of the plant.}
  \label{Fig666}
 \end{centering}
\end{figure*}
We summarize with the following theorem.
\begin{thm}\label{theorem3}
Consider the closed-loop system \eqref{ES_dynamicss}, and suppose that Assumptions \ref{assumptionFTGES} and \ref{assumption_plant} hold, and that the steady state input-to-output map \eqref{ssiom} satisfies Assumption \ref{assumption_1}. Then for any pair of admissible parameters $(q_1,q_2)$, $k>0$, $\delta>\nu>0$, there exists $\eta>0$ and $\varepsilon_2^*>0$ such that for all $\varepsilon_2\in(0,\varepsilon_2^*)$ there exists $a^*>0$ such that for all $a\in(0,a^*)$ there exists $\varepsilon_1^*>0$ such that for all $\varepsilon_1\in(0,\varepsilon_1^*)$ there exists $\varepsilon_0^*>0$ such that for all $\varepsilon_0\in(0,\varepsilon_0^*)$ system \eqref{ES_dynamicss} with initial conditions  $(u(0),\xi(0),\mu(0),x(0))\in (\{z^*\}+\delta\mathbb{B})\times \eta\mathbb{B}\times \mathbb{S}^n\times \Xi$ generates complete solutions and each trajectory satisfies $|z(\tau)-z^*|\leq \nu$, for all $\tau\in\text{dom}(u,\xi,\mu,x)$, such that $\tau\geq T_G^*$. \QEDB
\end{thm}
It is important to note that the fixed time convergence result of Theorem \ref{theorem3} holds in the $\tau$-time scale. Since $\tau=t\varepsilon_0$, the convergence result in the original time scale would translate to $t\geq T_G^*/\varepsilon_0$. Thus, smaller values of $\varepsilon_0$ generate larger values for the upper bound of the convergence time. This observation illustrates the key difference between fixed-time extremum seeking in dynamic plants versus static maps.
%
%

We finish this section by presenting the convergence result for the FTNES algorithm applied to the plant dynamics \eqref{plant_dynamics}. In this case, the closed-loop system is shown in Figure \ref{Fig66}, where the gains $k_i$ for $i\in\{1,2,3\}$ are defined as in the FTGES of \eqref{ES_dynamic}. Since the proof is almost identical to the proof of Theorem \ref{theorem3} by using the result of Theorem \ref{main_theorem2}, we present this result as a Corollary.
\begin{cor}\label{corollaryNewton}
Consider the closed-loop system shown in Figure \ref{Fig66} in the $\tau$-time scale, and suppose that Assumptions \ref{assumptionFTGES} and \ref{assumption_plant} hold and that the steady state input-to-output mapping \eqref{ssiom} satisfies Assumption \ref{assumption_2}. Then, for admissible parameters $(q_1,q_2)$, $k>0$, and each $\nu>0$ there exists $\eta>0$ and $\varepsilon_2^*>0$ such that for each $\varepsilon_2\in(0,\varepsilon_2^*)$ there exists $a^*>0$ such that for each $a\in(0,a^*)$ there exists $\varepsilon_1^*>0$ such that for each $\varepsilon_1\in(0,\varepsilon_1^*)$ there exists $\varepsilon_0^*$ such that for each $\varepsilon_0\in(0,\varepsilon_0^*)$ there exists a neighborhood $\mathcal{N}$ of $p^*:=(z^*, \nabla^2 \phi(z^*)^{-1},\nabla\phi(z^*),\ell(z^*))$ such that every solution with $(u(0),\xi(0),x(0),\mu)\in \mathcal{N}\times\mathbb{S}^n$ is defined for all time $\tau\geq0$, and satisfies  $|z(\tau)-z^*|\leq \nu$, for all $\tau\geq T_N^*$. \QEDB
\end{cor}
As in the gradient-based case, the convergence result of Corollary \ref{corollaryNewton} would imply that in the $t$-time scale the input satisfies $|z(t)-z^*|\leq \nu$ for all $t\geq T_N^*/\varepsilon_0$. Thus, as $\varepsilon_0\to0^+$ the convergence time grows unbounded since the controller is ``slowed down'' to guarantee enough time scale separation with respect to the plant dynamics \eqref{plant_dynamics}. Nevertheless, for the FTNES dynamics the value of $T^*_{N}$ does not depend on the cost function $\phi$.
\subsection{Numerical Example}
To illustrate the application of the fixed-time extremum seeking controller in dynamic plants, we consider the following dynamical system
\begin{equation}\label{dynamic_plant_example}
\begin{split}
\dot{x}_1&=-20x_1+5x_2+5z_1\\
\dot{x}_2&=-20x_2+5z_2,
\end{split}
\end{equation}
with output function given by
\begin{equation*}
y=10x_1^2+10x_2^2+\frac{x_1}{2}+\frac{x_2}{5}.
\end{equation*}
The quasi-steady state manifold is given by $\ell(z)=\left(\frac{z_2}{16}+\frac{z_1}{4}, \frac{z_2}{4}\right)$, and the steady state input-to-output mapping is 
\begin{equation*}
\phi(z)=h(\ell(z))=10\left(\frac{z_1}{4}+\frac{z_2}{16}\right)^2+\frac{5z_2^2}{8}+\frac{13z_2}{160}+\frac{z_1}{8},
\end{equation*}
which is strongly convex and has a minimum at $z^*=(-0.09,-0.04)$. We implement the FTGES using only measurements $y$ of the output plant, and parameters $\varepsilon_0=0.1$, $\varepsilon_1=0.0015$, $\varepsilon_2=0.05$, $q_1=3$, $q_2=1.5$, and $k=0.2$. We simulated 70 different trajectories of the closed-loop system, each trajectory with an initial condition $u(0)$ selected randomly from a compact set $\Omega_{10}$ using a uniform distribution. Figure \ref{Fig666} shows the time history of the resulting 70 trajectories, as well as the theoretical upper bound on the convergence time in the $\tau$-time scale normalized by $\varepsilon_0$. The insets show the evolution of the states $x_1$ and $x_2$, which were initialized at the origin, as well as the residual error after the convergence time.  As it can be observed, all the trajectories of $u$ converge to a small neighborhood of $z^*$ at a time that is upper bounded by the theoretical bound $T_G^*$ normalized by $\varepsilon_0$. Since the plant is stable, the states $x_1$ and $x_2$ eventually also converge to a  neighborhood of the quasi-steady state manifold $\ell(z^*)$.
\subsection{Connections with Existing Results in the Literature}
\label{subsec_connections}
We now discuss some connections between the ES dynamics presented in this paper and the existing gradient and Newton-based ES algorithms considered in \cite{DerivativesESC} and \cite{Newton}. In particular, we show that when $\alpha_1=\alpha_2=0$, the FTGES and FTNES dynamics recover the existing schemes in the literature which have only asymptotic (semi-global practical) convergence properties.
%
\subsubsection{Gradient-Based Scheme} The proposed FTGES \eqref{ES_dynamics1} can be seen as a generalizations of the existing gradient-based ES algorithms. In particular, when $\alpha_1=\alpha_2=0$, the FTGES dynamics applied to the dynamic plant \eqref{plant_dynamics} become
\begin{align}\label{ES_dynamicA}
\left(\begin{array}{c}
\dot{u}\vspace{0.1cm}\\
\dot{\xi}\vspace{0.1cm}\\
\dot{\mu}\vspace{0.1cm}\\
\dot{x}\vspace{0.1cm}
\end{array}\right)=\left(\begin{array}{c}
-2k_1\xi\vspace{0.1cm}\\
-k_2\big(\xi-F_G(y,\mu)\big)\vspace{0.1cm}\\
-k_3\mathcal{R}_{\kappa}\mu\vspace{0.1cm}\\
f(x,u+a\mathcal{D}u)
\end{array}\right).
\end{align}
In this case, by setting $k=0.5$, the reduced average dynamics of \eqref{ES_dynamicA} are given by
\begin{equation}\label{simplified_gradient_flow}
\dot{u}_r=-\nabla \phi(u_r),
\end{equation}
which is just a classic gradient descent. Under condition \eqref{condition2}, the function $\phi(u_r)$ is invex \cite{PLInequality}, which implies that $\mathcal{A}:=\{u^{*}_r\in\mathbb{R}^n: \nabla\phi(u^{*}_r)=0\}=\{u^{*}_r\in\mathbb{R}^n: u^{*}_r=\text{arg}\min_{u_r\in\mathbb{R}^n} \phi(u_r)\}$, i.e., every critical point is a global minimizer \cite{Invexity1}. Thus, the Lyapunov function $V=\phi(u_r)-\phi^*$ satisfies $\dot{V}=-|\nabla \phi(u_r)|^2$, which is zero only at points that minimize the function $\phi$. Since the cost function is radially unbounded and $\phi$ attains its minimum value, the set $\mathcal{A}$ is compact and UGAS under the dynamics \eqref{simplified_gradient_flow}. If one further assumes that $\nabla \phi$ is globally Lipschitz, the set $\mathcal{A}$ is indeed uniformly globally exponentially stable (UGES). Thus, in this case the $\mathcal{KL}$ bound that characterizes the convergence of $u$ in the original dynamics will be of the form $\beta(s_1,s_2)=c_1s_1\exp(-c_2s_2)$, for some $c_1,c_2>0$, see also \cite{zero_order_poveda_Lina}. Therefore, the choice $\alpha_1=\alpha_2=0$ recovers the existing results for gradient-based ES of \cite{tan06Auto} and \cite{DerivativesESC}. Finally, we note that when $\alpha_1=\alpha_2=1$, and $u_r$ is a scalar, we recover the finite-time gradient-based ES algorithm presented in \cite[Sec. 6.1]{Poveda:16}, which is discontinuous at the set of minimizers, and which has a convergence time that is dependent on the initial conditions of the optimizing state $u$.

\subsubsection{Newton-based Schemes} For the Newton-based case, setting $\alpha_1=\alpha_2=0$ in the FTNES dynamics applied to the plant \eqref{plant_dynamics} results in the closed-loop system:
\begin{align}\label{ES_dynamics1NB}
\left(\begin{array}{c}
\dot{u}\\\\
\dot{\xi}_1\\\\
\dot{\xi}_2\\\\
\dot{\mu}
\end{array}\right)=-\left(\begin{array}{c}
2k_1\xi_1\xi_2\vspace{0.1cm}\\
k_2\Big(\xi_1F_H(\phi,\mu)\xi_1-\xi_1\Big),\vspace{0.1cm}\\
k_2\Big(\xi_2-F_G(\phi,\mu)\Big)\vspace{0.1cm}\\
k_3\mathcal{R}_{\kappa}\mu\vspace{0.1cm}\\
f(x,u+a\mathcal{D}u)
\end{array}\right).
\end{align}
By choosing $k=0.5\gamma$, the reduced average dynamics of \eqref{ES_dynamics1NB} become
\begin{equation*}
\dot{u}_r=- \gamma\nabla^2\phi(u_r)^{-1} \nabla \phi(u_r),
\end{equation*}
which is a classic Newton-flow. Under Assumption \ref{assumption_2}, this dynamics render the unique minimizer $z^*$ locally exponentially stable with rate of convergence proportional to the tunable gain $\gamma$, which is consistent with the results from \cite{Newton}.

\section{Conclusions and Outlook}
\label{Sec_Conclusions}
In this paper we presented a novel class of extremum seeking controllers that achieve fixed-time convergence with a convergence time that is independent of the optimizing state. We considered both model-free gradient-based algorithms and model-free Newton-based algorithms for different classes of cost functions, and we showed that for the gradient-based algorithm the upper bound on the convergence time of the extremum seeking controller depends on the parameters of the cost. On the other hand, for the Newton-based dynamics the upper bound is independent of the cost function and can be prescribed a priori. Both results were established by using generalized averaging theory for non-smooth systems. We also extended these extremum seeking controllers to dynamical systems that generate well-defined steady state input-to-output mappings, and for which the fixed-time convergence property holds after a time scale transformation is performed. Our results were further validated by means of different single-variable and multivariable numerical examples.

The results of this paper open the door to novel opportunities for the development of other type of fixed-time extremum seeking controllers. In particular, the results of this paper can be extended and generalized for the solution of constrained extremum seeking problems, Nash seeking problems in games, tracking of time-varying optimizers, and hybrid extremum seeking controllers with fixed-time convergence properties. 

\vspace{0.5cm}
\section*{Acknowledgments}
The first author would like to thank Andy Teel for fruitful discussions on the semi-global practical fixed-time convergence properties for nonsmooth singularly perturbed systems in standard form.

\bibliographystyle{ieeetr}
\bibliography{Biblio}

\begin{thebibliography}{10}

\bibitem{ACCpovedaKrstic}
J.~I. Poveda and M.~Krsti\'{c}, ``An extremum seeking controller with
  semi-global practical prescribed finite-time convergence,'' {\em Submitted to
  the American Control Conference on September 26th, 2019, available online at:
  https://sites.google.com/site/jorgeivanpoveda/ftes}.

\bibitem{PovedaKrsticIFACWC20}
J.~I. Poveda and M.~Krsti\'{c}, ``Fixed-time newton-based extremum seeking,''
  {\em Submitted to 21st IFAC World Congress on Oct. 28th, 2019, available
  online at: https://sites.google.com/site/jorgeivanpoveda/ftes}.

\bibitem{fixed_time}
K.~Garg and D.~Panagou, ``Fixed-time stable gradient-flow schemes: Applications
  to continuous-time optimization,'' {\em arXiv:1808.10474}, 2018.

\bibitem{KrsticBookESC}
K.~B. Ariyur and M.~Krsti\'{c}, {\em Real-Time Optimization by Extremum-Seeking
  Control}.
\newblock Wiley, 2003.

\bibitem{DerivativesESC}
D.~Ne\u{s}i\'{c}, Y.~Tan, W.~H. Moase, and C.~Manzie, ``A unifying approach to
  extremum seeking: Adaptive schemes based on estimation of derivatives,'' {\em
  49th IEEE Conference on Decision and Control}, pp.~4625--4630, 2010.

\bibitem{Moase}
W.~Moase, Y.~Tan, D.~Ne\u{s}i\'c, and C.~Manzie, ``Non-local stability of a
  multi-variable extremum-seeking scheme,'' {\em In Proc. of IEEE Australian
  Control Conference}, pp.~38--43, 2011.

\bibitem{Durr:Lie}
H.~Durr, M.~S. Stankovic, C.~Ebenbauer, and K.~H. Johansson, ``Lie bracket
  approximation of extremum seeking systems,'' {\em Automatica}, vol.~49,
  pp.~1538--1552, 2013.

\bibitem{tan06Auto}
Y.~Tan, D.~Ne\v{s}i\'{c}, and I.~M. Mareels, ``On non-local stability
  properties of extremum seeking control,'' {\em Automatica}, vol.~42, no.~6,
  pp.~889--903, 2006.

\bibitem{GrushkovskayaLie}
V.~Grushkovskaya, A.~Zuyev, and C.~Ebenbauer, ``On a class of generating vector
  fields for the extremum seeking problem: Lie bracket approximation and
  stability properties,'' {\em Automatica}, vol.~94, pp.~151--160, 2018.

\bibitem{PowerES}
A.~Ghaffari, M.~Krsti\'{c}, and S.~Seshagiri, ``Power optimization for
  photovoltaic micro-converters using multivariable gradient-based
  extremum-seeking,'' {\em IEEE Transactions on Control System Technology},
  vol.~22, no.~6, pp.~2141--2149, 2014.

\bibitem{Poveda:15}
J.~I. Poveda and N.~Quijano, ``Shahshahani gradient-like extremum seeking,''
  {\em Automatica}, vol.~58, pp.~51--59, 2015.

\bibitem{HeavyBollES}
S.~Michalowsky and C.~Ebenbauer, ``The multidimensional n-th order heavy ball
  method and its application to extremum seeking,'' {\em 53rd IEEE Conf.
  Decision Control}, pp.~2660--2666, 2014.

\bibitem{zero_order_poveda_Lina}
J.~I. Poveda and N.~Li, ``Robust hybrid zero-order optimization algorithms with
  acceleration via averaging in continuous time,'' {\em arXiv:1909.00265},
  2019.

\bibitem{Strizic:17_CDC}
T.~Strizic, J.~I. Poveda, and A.~R. Teel, ``Hybrid gradient descent for robust
  global optimization on the circle,'' {\em IEEE 56th Conference on Decision
  and Control}, pp.~2985--2990, 2017.

\bibitem{ThiagoKrstic}
T.~Oliveira, M.~Krsti\'{c}, and D.~Tsubakino, ``Extremum seeking for static
  maps with delays,'' {\em IEEE Trans. Autom. Control}, vol.~62, no.~4,
  pp.~1911--1926, 2017.

\bibitem{Grushkovskaya2017}
V.~Grushkovskaya, H.~Durr, C.~Ebenbauer, and A.~Zuyev, ``Extremum seeking for
  time-varying functions using lie bracket approximations,'' {\em 20th {W}orld
  {C}ongress}, vol.~50, no.~1, pp.~522--5528, 2017.

\bibitem{Suttner2017}
R.~Suttner and S.~Dashkovskiy, ``Exponential stability for extremum control
  systems,'' {\em 20th {IFAC} {W}orld {C}ongress}, pp.~15464--15470, 2017.

\bibitem{Poveda:16}
J.~I. Poveda and A.~R. Teel, ``A framework for a class of hybrid extremum
  seeking controllers with dynamic inclusions,'' {\em Automatica}, vol.~76,
  pp.~113--126, 2017.

\bibitem{Guay:03}
M.~Guay and T.~Zhang, ``Adaptive extremum seeking control of nonlinear dynamic
  systems with parametric uncertainties,'' {\em Automatica}, vol.~39,
  pp.~1283--1293, 2003.

\bibitem{Guay:15}
M.~Guay and D.~Dochain, ``A time-varying extremum-seeking control approach,''
  {\em Automatica}, vol.~51, pp.~356--363, 2015.

\bibitem{AttaGuay18}
K.~T. Atta, A.~Johansson, and M.~Guay, ``On the generalization and stability
  analysis of pareto seeking control,'' {\em IEEE Control Systems Letters},
  vol.~2, no.~1, pp.~145--150, 2018.

\bibitem{PhasorGuay}
K.~T. Atta and M.~Guay, ``Adaptive amplitude fast proportional integral phasor
  extremum seeking control for a class of nonlinear system,'' {\em Journal of
  Process Control}, vol.~83, pp.~147--154, 2019.

\bibitem{improvedES}
K.~T. Atta, A.~Johansson, and T.~Gustafsson, ``Accuracy improvement of extremum
  seeking control,'' {\em IEEE Transactions on Automatic and Control}, vol.~62,
  no.~4, pp.~1952--1958, 2017.

\bibitem{Fish_ESC}
J.~Cochran, E.~Kanso, S.~D. Kelly, H.~Xiong, and M.~Krsti\'c, ``Source seeking
  for two nonholonomic models of fish locomotion,'' {\em IEEE Trans. on
  Robotics}, vol.~25, pp.~1166--1176, 2009.

\bibitem{Frihauf:12}
P.~Frihauf, M.~Krsti\'{c}, and T.~Basar, ``Nash equilibrium seeking in
  noncooperative games,'' {\em IEEE Transactions on Automatic Control},
  vol.~57, no.~5, pp.~1192--1207, 2012.

\bibitem{PovedaCDC17_a}
J.~I. Poveda, P.~N. Brown, J.~R. Marden, and A.~R. Teel, ``A class of
  distributed adaptive pricing mechanisms for societal systems with limited
  information,'' {\em 56th IEEE Conference on Decision and Control},
  pp.~1490--1495, 2017.

\bibitem{Kutadinata:14_Traffic}
R.~Kutadinata, W.~Moase, C.~Manzie, L.~Zhang, and T.~Garoni, ``Enhancing the
  performance of existing urban traffic light control through
  extremum-seeking,'' {\em Transportation Research Part C: Emerging
  Technologies}, vol.~62, pp.~1--20, 2016.

\bibitem{NonC2Krstic}
A.~Scheinker and M.~Krsti\'{c}, ``Non-{C}2 lie bracket averaging for nonsmooth
  extremum seekers,'' {\em Journal of Dynamic Systems, Measurement, and
  Control}, vol.~136, no.~1, pp.~1--10, 2014.

\bibitem{JaimeFixed_Time}
E.~Cruz-Zavala, J.~A. Moreno, and L.~Fridman, ``Uniform second-order sliding
  mode observer for mechanical systems.,'' {\em Proc. Int. Workshop Variable
  Struct. Syst}, pp.~14--19, 2010.

\bibitem{finite_timeEngelTAC}
R.~Engel and G.~Kreisselmeier, ``A continuous-time observer which converges in
  finite time,'' {\em IEEE Transactions on Automatic and Control}, vol.~47,
  pp.~1202--1204, 2002.

\bibitem{Fixed_timeTAC}
A.~Polyakov, ``Nonlinear feedback design for fixed-time stabilization of linear
  control systems,'' {\em IEEE Transactions on Automatic and Control}, vol.~57,
  no.~8, pp.~2106--2110, 2012.

\bibitem{OldFinite_Time}
M.~James, ``Finite time observers and observability,'' {\em in Proc. of
  Conference on Decision and Control}, pp.~770--771, 1990.

\bibitem{PralyFinite_Time}
V.~Andrieu, L.~Praly, and A.~Astolfi, ``Homogeneous approximation, recursive
  observer design, and output feedback,'' {\em SIAM J. Control Optim}, vol.~47,
  no.~4, pp.~1814--1850, 2008.

\bibitem{romeronips}
O.~Romero and M.~Benosman, {\em Finite-Time Convergence of Continuous-Time
  Optimization Algorithms via Differential Inclusions}.
\newblock Neurips, to appear., 2019.

\bibitem{Garg_Inequalities}
K.~Garg, M.~Baranwal, R.~Gupta, R.~Vasudevan, and D.~Panagou, ``Fixed-time
  stable proximal dynamical system for solving mixed variational inequality
  problems,'' {\em arXiv:1908.03517}, 2019.

\bibitem{Wang:12_Automatica}
W.~Wang, A.~Teel, and D.~Ne\u{s}i\'{c}, ``Analysis for a class of singularly
  perturbed hybrid systems via averaging,'' {\em Automatica}, vol.~48, no.~6,
  pp.~1057--1068, 2012.

\bibitem{averaging_singularHDS}
W.~Wang, A.~R. Teel, and D.~Ne\u{s}\'ic, ``Averaging in singularly perturbed
  hybrid systems with hybrid boundary layer systems,'' {\em 51st IEEE
  Conference on Decision and Control}, pp.~6855--6860, 2012.

\bibitem{Newton}
A.~Ghaffari, M.~Krsti\'{c}, and D.~Ne{\v{s}}i{\'c}, ``Multivariable
  newton-based extremum seeking,'' {\em Automatica}, vol.~48, pp.~1759--1767,
  2012.

\bibitem{Beck2014}
A.~Beck, {\em Introduction to Nonlinear Optimization}.
\newblock MOS-SIAM Series on Optimization, 2014.

\bibitem{HDS}
R.~Goebel, R.~Sanfelice, and A.~R. Teel, {\em Hybrid Dynamical System}.
\newblock Princeton, NJ: Princeton University Press, 2012.

\bibitem{khalil}
H.~K. Khalil, {\em Nonlinear Systems}.
\newblock Upper Saddle River, NJ: Prentice Hall, 2002.

\bibitem{durrEbenbauer}
H.~D\"{u}rr, M.~Krsti\'{c}, A.~Scheinker, and C.~Ebenbauer, ``Extremum seeking
  for dynamic maps using liebrackets and singular perturbations,'' {\em
  Automatica}, vol.~83, no.~91-99, 2017.

\bibitem{Krstic2000}
M.~Krsti\'{c} and H.~Wang, ``Stability of extremum seeking feedback for general
  nonlinear dynamic systems,'' {\em Automatica}, vol.~36, no.~4, pp.~595--601,
  2000.

\bibitem{TanAndNesic2006Local}
Y.~Tan, D.~Ne\u{s}i\'c, and I.~Mareels, ``On non-local stability properties of
  extremum seeking controllers,'' {\em Automatica}, vol.~42, no.~6,
  pp.~889--903, 2006.

\bibitem{NewtonESC}
C.~Labar, E.~Garone, M.~Kinnaert, and C.~Ebenbauer, ``Newton-based extremum
  seeking: A second-order lie bracket approximation approach,'' {\em
  Automatica}, vol.~105, pp.~356--367, 2019.

\bibitem{NewtonSemiglobal}
F.~Alvarez, H.~Attouch, J.~Bolte, and P.~Redont, ``A second-order gradient-like
  dissipative dynamical system with hessian-driven damping: {A}pplication to
  optimization and mechanics,'' {\em J. Math. Pures Appl.}, vol.~81,
  pp.~747--779, 2002.

\bibitem{black_box_CDC}
J.~I. Poveda, R.~Kutadinata, C.~Manzie, D.~Ne\u{s}i\'{c}, A.~R. Teel, and C.~K.
  Liao, ``Hybrid extremum seeking for black-box optimization in hybrid plants:
  An analytical framework,'' {\em IEEE Conference on Decision and Control},
  pp.~2235--2240, 2018.

\bibitem{PLInequality}
H.~Karimi, J.~Nutini, and M.~Schmidt, {\em Linear Convergence of Gradient and
  Proximal-Gradient Methods Under the Polyak-Lojasiewicz Condition}, vol.~9851,
  ch.~Machine Learning and Knowledge Discovery in Databases, ECML PKDD 2016.
\newblock Springer, 2016.

\bibitem{Invexity1}
A.~Ben-Israel and B.~Mond, ``What is invexity?,'' {\em J. Austral. Math. Soc.
  (Series B)}, vol.~28, pp.~1--9, 1986.

\end{thebibliography}

\end{document}